\documentclass[onefignum,onetabnum]{siamart220329}

\usepackage{subfigure,threeparttable}
\usepackage{booktabs,multirow}
\usepackage{lipsum,enumerate,enumitem}
\usepackage{amsfonts,mathrsfs}
\usepackage{graphicx}
\usepackage{cite}
\usepackage{epstopdf}
\usepackage{algorithmic}
\ifpdf
  \DeclareGraphicsExtensions{.eps,.pdf,.png,.jpg}
\else
  \DeclareGraphicsExtensions{.eps}
\fi

\newsiamremark{remark}{Remark}
\newsiamremark{hypothesis}{Hypothesis}
\crefname{hypothesis}{Hypothesis}{Hypotheses}
\newsiamthm{claim}{Claim}

\def\es{\mathrm{es}}
\def\ck{{\scriptstyle{\mathcal{C}}}}
\newcommand{\cte}{{\mathrm{T}}}
\newcommand{\cnt}{{\mathrm{N}}}

\newcommand{\zd}{\,\mathrm{d}}

 \newcommand{\abst}[1]{|#1|}
 \newcommand{\bra}[1]{\left(#1\right)}
\newcommand{\brab}[1]{\big(#1\big)}
\newcommand{\braB}[1]{\Big(#1\Big)}
\newcommand{\brat}[1]{(#1)}
\newcommand{\kbra}[1]{\left[#1\right]}
\newcommand{\kbrab}[1]{\big[#1\big]}
\newcommand{\kbraB}[1]{\Big[#1\Big]}

\newcommand{\myinnert}[1]{\langle#1\rangle}
\newcommand{\myinnerb}[1]{\big\langle#1\big\rangle}
\newcommand{\myinnerB}[1]{\Big\langle#1\Big\rangle}
 \newcommand{\mynorm}[1]{\left\|#1\right\|}
 \newcommand{\mynormb}[1]{\big\|#1\big\|}
 
 \newcommand{\mynormt}[1]{\|#1\|}

\def\lan#1{\textcolor{blue}{#1}}

\headers{Original energy laws of IERK, EERK and CIFRK}{X. Wang, X. Zhao and H.-L. Liao}

\title{A unified framework on the original energy laws of three effective classes of Runge-Kutta methods for phase field crystal type models\thanks{Submitted to the editors \today.
		\funding{This work is supported by the NSF of China under grants 12471383, 12071216 and U22B2046, Jiangsu Provincial Scientific Research Center of Applied Mathematics under grant BK20233002, and the Natural Science Foundation of Jiangsu Province under grant BK20222003.}}}

\author{Xuping Wang\thanks{School of Mathematics, Nanjing University of Aeronautics and Astronautics, Nanjing 211106, China. (\email{wangxp@nuaa.edu.cn})}
	\and Xuan Zhao\thanks{School of Mathematics, Southeast University, Nanjing 210096, China. (\email{xuanzhao11@seu.edu.cn})}\and Hong-lin Liao\thanks{Corresponding author. ORCID 0000-0003-0777-6832. School of Mathematics, Nanjing University of Aeronautics and Astronautics, Nanjing 211106, China;
		Key Laboratory of Mathematical Modeling and High Performance Computing of Air Vehicles (NUAA), MIIT, Nanjing 211106, China. (\email{liaohl@nuaa.edu.cn}, \email{liaohl@csrc.ac.cn})} }

\usepackage{amsopn}

\begin{document}
  
\maketitle
  
\begin{abstract}
The main theoretical obstacle to establish the original energy dissipation laws of Runge-Kutta methods for phase-field equations is to verify the maximum norm boundedness of the stage solutions without assuming global Lipschitz continuity of the nonlinear bulk. We present a unified theoretical framework for the energy stability of three effective classes of Runge-Kutta methods, including the additive implicit-explicit Runge-Kutta, explicit exponential Runge-Kutta and corrected integrating factor Runge-Kutta methods, for the Swift-Hohenberg and phase field crystal models. By the standard discrete energy argument, it is proven that the three classes of Runge-Kutta methods preserve the original energy dissipation \lan{laws} if the associated differentiation matrices are positive definite. Our main tools include the differential form with the associated differentiation matrix, the discrete orthogonal convolution \lan{kernels} and the principle of mathematical induction. Many existing Runge-Kutta methods in the literature are revisited by evaluating the lower bound on the minimum eigenvalues of the associated differentiation matrices. Our theoretical approach paves a new way for the internal nonlinear stability of  Runge-Kutta methods for dissipative semilinear parabolic problems.
\end{abstract}

  
\section{Introduction}
\setcounter{equation}{0}

We suggest a unified theoretical framework to establish the original  energy dissipation laws at all stages for three computationally effective classes of  Runge-Kutta (RK) methods, including  the additive implicit-explicit RK (IERK), explicit exponential RK (EERK) or exponential time differencing RK (ETDRK) and corrected integrating factor RK (CIFRK) methods, in solving the phase field crystal (PFC) type model 
\begin{equation}\label{model: PFC type}
  \partial_t \Phi= \mathcal{M} \kbra{(I+\Delta)^2 \Phi - f(\Phi)}.
  \end{equation}
It is a very powerful model to describe crystal dynamics at the atomic scale in space and on diffusive scales in time, especially for crystalline self-organization or pattern formation on the atomistic scale, cf. \cite{ElderGrant:2004SPFC,PFC:2002PRL,ProvatasElder:2011,SH:1977PRA}.
Here, $\Phi$ represents the periodic density field in the domain $\Omega=(0, L)^2 \subset \mathbb{R}^2$, $f(\Phi) := \varepsilon \Phi-\Phi^3$ is the nonlinear bulk with a constant $\varepsilon \in (0, 1)$ related to the temperature, and the mobility operator $\mathcal{M}$ is negative definite. The PFC type model \eqref{model: PFC type} can be formulated as the gradient flow $\partial_t \Phi = \mathcal{M} \frac{\delta E}{\delta \Phi}$
with the Swift-Hohenberg free energy functional \cite{SH:1977PRA}
\begin{equation}\label{eq: energy functional}
  E[\Phi]:= \int_{\Omega} \kbraB{\tfrac12 \Phi(I+\Delta)^2\Phi + F(\Phi)} \zd \mathbf{x},
\end{equation}
in which the double-well potential $F(\Phi) = \frac14 \Phi^4-\frac\varepsilon2 \Phi^2$ satisfies $f(\Phi)=-F^\prime(\Phi)$. In this work, we consider two well-known cases:
\begin{itemize}[leftmargin=10mm]
  \item The $L^2$ gradient flow with $\mathcal{M} = -I$ in \eqref{model: PFC type} arrives at  the Swift-Hohenberg (SH) equation \cite{SH:1977PRA}, that is, $\partial_t\Phi= - (I+\Delta)^2 \Phi +f(\Phi)$.
  \item The  $H^{-1}$ gradient flow with $\mathcal{M} = \Delta$ in \eqref{model: PFC type} gives the volume-conservative PFC equation \cite{PFC:2002PRL}, that is, $\partial_t\Phi = \Delta \kbrab{(I+\Delta)^2 \Phi - f(\Phi)}$.
\end{itemize}
As is well known, the long-time dynamics approaching the steady-state solution $\Phi^*$, that is, $(I+\Delta)^2 \Phi^* = f(\Phi^*)$, of the dissipative system \eqref{model: PFC type} satisfies the following original energy dissipation law
\begin{equation}\label{eq: energy dissipation law}
  \frac{\mathrm{d} E}{\mathrm{d} t}=\braB{\frac{\delta E}{\delta \Phi}, \partial_t\Phi}=(\mathcal{M}^{-1}\partial_t\Phi, \partial_t\Phi)\le0,
\end{equation}
where $(f, g):=\int_{\Omega} f g \mathrm{d} \mathbf{x}$ with the $L^{2}$ norm $\mynorm{f}_{L^2}:=\sqrt{(f, f)}$ for all $f, g \in L^{2}(\Omega)$.

In this paper, we will consider some linearly energy stable RK methods with the explicit approximation of the nonlinear bulk $f$, which are always computationally efficient by removing the inner iteration at each time stage.  The linear stabilization technique \cite{DuJuLiQiao:2021SIREV,FuTangYang:2024MCOM} is adopted by 
introducing the following stabilized operators,
\begin{align}\label{def: stabilized parameter}
	\mathcal{L}_{\kappa}\Phi:=(I+\Delta)^2 \Phi+\kappa \Phi \quad\text{and}\quad f_{\kappa}(\Phi):=\kappa \Phi+f(\Phi),
\end{align} 
where $\kappa\ge0$ represents the stabilized parameter controlling the possible instability from the explicit approximation of $f$.
Then the PFC type problem \eqref{model: PFC type} becomes 
\begin{align}\label{model: PFC type stabilized}
	\partial_t \Phi=\mathcal{M} \kbra{\mathcal{L}_{\kappa}\Phi-f_{\kappa}(\Phi)}\quad\text{for $\lan{\mathbf{x}}\in\Omega$},
\end{align}
which is equivalent to \eqref{model: PFC type} at the continuous level.

In the existing literature, there are many numerical methods developed for the PFC type model \eqref{model: PFC type}. As pointed out in \cite{LiQiaoWangZheng:2024arxiv}, a theoretical justification of the energy dissipation law \eqref{eq: energy dissipation law} has always been recognized as an important mathematical indicator for an effective numerical scheme to the gradient flow \lan{system}, since the energy stability always plays a crucial role in the long-time numerical simulations. Actually, for the PFC type equation \eqref{model: PFC type}, there have been extensive works to develop energy stable time-stepping methods, cf. \cite{ChengWangWise:2019CCP,CuiNiuXu:2024JSC,Lee:2019CMAME,LiaoJiZhang:2022IMA,LiaoKang:2023IMA,LiuYin:2019JSC,SuFangetal:2019CAM,SunZhaoetal:2022CNSNS,SunZhangetal:2024arxiv,WiseWangetal:2009SINUM} and the references therein, with certain modified discrete energy functional by adding some nonnegative small terms to the original energy \eqref{eq: energy functional}, especially for high-order multi-step approaches. 

On the other hand, linear time-stepping methods preserving the decay of original energy $E[\Phi]$ would be more attractive in both theoretical and practical manner. Very recently, there \lan{is} some key progress in the theoretical justification of original energy dissipation for several computationally efficient classes of RK methods in the long-time simulation of gradient flow problems. Maybe the first justification was done by Fu and Yang \cite{FuYang:2022JCP}, in which the second-order ETDRK method in \cite{CoxMatthews:2002JCP} was shown to preserve the decay of original energy. 
By calculating the energy difference between adjacent time levels, a sufficient condition was derived in \cite{FuShenYang:2024SCM} to ensure the monotonic decay of original energy for a general class of ETDRK or EERK methods. The same sufficient condition was also imposed in \cite{LiaoWang:2024MCOM} from the perspective of differential form with the associated differentiation matrix of general EERK methods. However, the above studies \cite{FuYang:2022JCP,FuShenYang:2024SCM,LiaoWang:2024MCOM} are still theoretically limited since they always assume that the nonlinear bulk $f$ is globally Lipschitz continuous.

For specific gradient flow problems, it is always theoretically difficult to remove the priori assumption of global Lipschitz continuity for the original energy decay of multi-stage RK methods, cf. \cite{DongLiQiaoZhang:2023ANM,FuShenYang:2024SCM,FuTangYang:2024MCOM,FuYang:2022JCP}, because it requires the uniform boundedness of discrete solutions at all stages (establishing the maximum bound principle of stage solutions is an alternative strategy \cite{DuJuLiQiao:2019SINUM,DuJuLiQiao:2021SIREV} for the Allen-Cahn type phase field models). 
Benefiting from the existence of the sixth-order dissipation term, Li and Qiao \cite{LiQiao:2024SISC} verified the monotonic decay of original energy for the ETD1 and ETDRK2 methods in solving the PFC equation by establishing the $H^2$ norm error estimates and the uniform boundedness of solution with the help of embedding inequality. 
Very recently, Li et al. \cite{LiQiaoWangZheng:2024arxiv} suggested a novel technique, called global-in-time energy estimate, to the numerical analysis of the ETDRK2 scheme, cf. the case $c_2=1$ of \eqref{scheme: EERK2}, for the PFC equation and established the original energy dissipation property for any final time (the uniform boundedness of solution is then also valid for any final time). The basic approach of the global-in-time energy estimate is: (Step 1) make priori assumption of decreasing energy at the previous time step, so that the discrete $H^2$ and $L^\infty$ bounds of the numerical solution become available at the current step; (Step 2) by reformulating each RK stage into two sub-stages and applying some delicate eigenvalue estimates of nonlinear terms in the Fourier pseudo-spectral space, the numerical system at the intermediate time stage and the next time step is carefully analyzed for a single-step $H^2$ estimate, which recovers the $L^\infty$ bounds of the numerical solution at the intermediate stage and the next time step. \lan{After the submission of this article, we are informed that the technique of global-in-time energy estimate has been applied by Zhang et al. \cite{ZhangWangTeng:2024SINUM} to a second-order Radau-type IERK method \cite[Section 3.6]{AscherRuuthSpiteri:1997} for the PFC equation. Sun et al. \cite{SunZhangetal:2024arxiv} used almost the same technique to examine a weak variant of ETDRK2 method, cf. the case $c_2=\tfrac{1}{2}$ of \eqref{scheme: EERK2-w}, for the SH equation.} 
Nonetheless, this technical complexity makes it difficult to generalize the global-in-time energy estimate to some high-order or parameterized EERK methods in \cite{HochbruckOstermann:2005SINUM,HochbruckOstermann:2010EERK}.
 
As a closely related approach to the above EERK or ETDRK methods, the so-called integrating factor RK (IFRK) or Lawson methods \cite{Lawson:1967SINUM} were also applied in {recent} years to the Allen-Cahn type gradient flows since they have obvious advantage in maintaining the maximum bound principle, see  \cite{JuLiQiaoYang:2021JCP,LiLiJuFeng:2021SISC,ZhangYanQianSong:2022CMAME}. 
However, IFRK methods have the so-called exponential effect (the maximum bound of solution decays exponentially when the time-step  size is properly large) and always can not preserve the original energy dissipation law \eqref{eq: energy dissipation law} although it may admit a modified energy law \cite[Theorem 2.1]{LiaoWangWen:2024JCP}. 
Actually, the method preserving the modified energy dissipation may not necessarily preserve the original energy dissipation law \cite{FuTangYang:2024MCOM,SAV:2022JCP}. {In recent work}  \cite{LiaoWangWen:2024JCP}, the exponential effect and the non-preservation of original energy dissipation law of IFRK methods were attributed to the non-preservation of steady-state, cf. \cite{CoxMatthews:2002JCP,Fornberg:1999JCP,HuangShu:2018JCP}. 
Some CIFRK methods were then developed in \cite{LiaoWangWen:2024JCP} by using two classes of difference correction, including the telescopic correction and nonlinear-term translation corrections, enforcing the preservation of steady-state solution $\Phi^*$ for any time-step sizes.
 Also, most of CIFRK methods in \cite{LiaoWangWen:2024JCP} were shown to preserve  the original energy dissipation law \eqref{eq: energy dissipation law} by using the associated differential form if the nonlinear bulk $f$ is globally Lipschitz continuous. Obviously, proving the uniform boundedness of discrete solutions at all stages generated by the CIFRK methods would be useful for specific gradient flow problems without the  priori assumption of global Lipschitz continuity.
 
It is to note that, the computational efficiency of the above two exponential classes of RK methods heavily relies on the computation of matrix exponentials. 
However, as pointed out in \cite{FasiGaudreaultLundSchweitzer:2024}, efficient algorithm to accurately compute the matrix exponential is still limited especially for general applications. Currently, the IERK methods have their own advantages in practical numerical simulations due to the lower computational complexity compared with these exponential classes of RK methods. By extensive numerical experiments in \cite{ShinLeeLee:2017CMA}, the IERK methods for gradient flow problems were observed to maintain the decay of original energy. This phenomenon was subsequently verified in \cite{ShinLeeLee:2017JCP} for a special class of nonlinear convex splitting RK methods. For the Radau-type (also known as ARS-type) IERK methods \cite{AscherRuuthSpiteri:1997,IzzoJackiewicz:2017}, cf. the form \eqref{IERK: general form} with the vanished coefficient $a_{i,1}=0$, Fu et al. \cite{FuShenYang:2024SCM} derived some sufficient conditions to maintain the decay of the original energy by computing the original energy difference between adjacent time levels. Very recently, by reformulating the IERK method into a steady-state preserving form, a unified theoretical framework was suggested in \cite{LiaoWangWen:2024IERK} to establish the original energy dissipation law of general IERK methods, including the Lobatto-type and Radau-type IERK methods. Nonetheless, the recent studies \cite{FuShenYang:2024SCM,LiaoWangWen:2024IERK} are also theoretically limited because they always assume that the nonlinear bulk $f$ is globally Lipschitz continuous. 
  
\textit{In summary, the latest developments on the IERK, EERK and CIFRK  methods for gradient flow problems have raised an urgent issue, that is, to build certain unified theory to establish uniform boundedness of solutions at all stages.} This task would be ambiguous at the first glance due to that the three methods have three different formulations in general. 
Fortunately, in recent studies \cite{LiaoWang:2024MCOM,LiaoWangWen:2024IERK,LiaoWangWen:2024JCP}, the above three effective classes of RK methods can be formulated into a unified differential form with three different types of lower triangle differentiation matrices $D(z)=\kbrab{d_{ij}(z)}_{i,j=1}^s,$ where $d_{ij}$ vanishes when $i<j$.
For simplicity,  the spatial operators $\Delta$, $\mathcal{L}_{\kappa}$ and $\mathcal{M}$ are approximated by the Fourier pseudo-spectral method, as described in \cref{section: theoretical framework}, with the associated discrete operators (matrices) $\Delta_h$, $L_{\kappa}$ and $M_{h}$.  For the PFC type model \eqref{model: PFC type stabilized}, the unified differential form reads 
\begin{align}\label{eq: unified differential form}
	\frac1{\tau}\sum_{j=1}^i d_{ij}(\tau M_h L_{\kappa}) \delta_{\tau} u_h^{n, j+1} = M_h \kbraB{L_{\kappa} u_h^{n, i+\frac12} - f_{\kappa}(u_h^{n, i})} \\
  = M_h \kbraB{(I+\Delta_h)^2 u_h^{n, i+\frac12} - f(u_h^{n, i}) + \tfrac\kappa2 \delta_{\tau} u_h^{n, i+1}} \notag
\end{align}
for $n\geq1$ and $1\le i\le s$, where the stage difference $\delta_{\tau} u_h^{n, j+1}:=u_h^{n, j+1}-u_h^{n, j}$, and the stage averaged operator $u_h^{n, \lan{j}+\frac12}:=(u_h^{n, j+1}+u_h^{n, j})/2$.  
For a given $s$-stage RK method, let $u_h^{n,i}$ be the numerical approximation of the exact solution $\Phi_h^{n,i}:=\Phi(\mathbf{x}_h, t_{n}+c_i\tau_n)$ at the abscissas $c_1:=0$, $c_i>0$ for $2 \leq i \leq s$, and $c_{s+1}:=1$. Let $\phi_h^n$ be the  approximation of $\Phi_h^n:=\Phi(\mathbf{x}_h, t_n)$ at the discrete time level $t_n$ for $n \geq 0$ such that the stage computation procedure \eqref{eq: unified differential form} starts from the previous level solution $u_h^{n,1}:=\phi_h^{n-1}$ and generates the next level solution by $\phi_h^{n}:=u_h^{n,s+1}$.
 
Always, a non-symmetric matrix $D$ is said to be positive (semi-)definite if the symmetric part $\mathcal{S}\brat{D}:=\brat{D+D^T}/2$ is positive (semi-)definite. Our main result is stated as follows. 
\begin{theorem}\label{theorem: energy stability}
  Assume that the nonlinear bulk $f$ is continuously differentiable and the differentiation matrix $D(z)$ is positive definite for any $z < 0$. There \lan{exist} two positive constants $\tau_0$ and $\ck_f$ such that when the time-step size $\tau\le \tau_0$ and the stabilized parameter $\kappa \geq \ck_f$, the $s$-stage RK methods having the general form \eqref{eq: unified differential form} for the PFC type model \eqref{model: PFC type stabilized} preserve the original energy dissipation law \eqref{eq: energy dissipation law}, 
  \begin{equation}\label{theoremResult: discrete energy law}
    E[u^{n,k+1}]-E[u^{n,1}]\le\frac1{\tau}\sum_{i=1}^k \myinnerB{M_h^{-1}\delta_{\tau} u^{n, i+1},\sum_{j=1}^i d_{ij}(\tau M_h L_{\kappa})\delta_{\tau} u^{n, j+1}}
  \end{equation}
  for $n\ge1$ and $1\le k\le s$. Moreover, it holds that
  \[
  \mynormb{(I+\Delta_h) u^{n,i}} + \mynormb{u^{n,i}} \leq \ck_0:=\sqrt{4E[\phi^{0}]+(1+\varepsilon)^2}\quad\text{for $n\ge1$ and $1\le i\le s$.}
  \]
\end{theorem}

In \cref{section: theoretical framework}, we describe the spatial discretization briefly and present the detail proof of Theorem \ref{theorem: energy stability} by the traditional discrete energy arguments including the $L^2$ norm analysis and the $H^2$ semi-norm analysis, in which we always assume that the minimum eigenvalue of the symmetric matrix $\mathcal{S}\kbra{D(z)}$ has a finite, $z$-independent lower bound $\lambda_{\es}>0$. Then, in \cref{section: lambda}, we will check the existence of such lower bounds $\lambda_{\es}$ by revisiting some specific RK methods, including the IERK, EERK and  \lan{CIFRK} methods. The last section presents some concluding remarks and further issues including the limitations of our theoretical framework.

\section{A concise framework for energy stability}\label{section: theoretical framework}

\subsection{Spatial approximation}
At first, we recall the Fourier pseudo-spectral method for spatial approximation. Consider the uniform length $h_x=h_y=h:=L/M$ in spatial directions for an even positive integer $M$. Define the discrete grid $\Omega_h:=\{\mathbf{x}_h=(mh,nh)\mid 1\leq m,n\leq M\}$ and $\bar{\Omega}_h:=\{\mathbf{x}_h=(mh,nh)\mid 0\leq m,n\leq M\}$. For a periodic function $v(\mathbf{x})$ on $\bar{\Omega}$, let $P_{M}: L^{2}(\Omega) \rightarrow \mathscr{F}_{M}$ be the standard $L^{2}$ projection operator onto the trigonometric polynomials space $\mathscr{F}_{M}$ (all trigonometric polynomials of degree up to $M / 2$), and $I_{M}: L^{2}(\Omega) \rightarrow \mathscr{F}_{M}$ be the trigonometric interpolation operator \cite{ShenTangWang:2011book}, that is,
\[
  (P_{M} v)(\mathbf{x})=\sum_{m, n=-M / 2}^{M / 2-1} \widehat{v}_{m, n} e_{m, n}(\mathbf{x}), \quad (I_{M} v)(\mathbf{x})=\sum_{m, n=-M / 2}^{M / 2-1} \tilde{v}_{m, n} e_{m, n}(\mathbf{x}),
\]
where the complex exponential basis functions $e_{m, n}(\mathbf{x}):=e^{\mathrm{i} \nu(m x+n y)}$ with $\nu=2 \pi / L$. The coefficients $\widehat{v}_{m, n}$ refer to the standard Fourier coefficients of $v(\mathbf{x})$, and the pseudo-spectral coefficients $\widetilde{v}_{m, n}$ are determined such that $(I_{M} v)(\mathbf{x}_{h})=v_{h}$. The Fourier pseudo-spectral first and second order derivatives of $v_{h}$ are given by
\[
  \mathcal{D}_{x} v_{h}:=\sum_{m, n=-M / 2}^{M / 2-1}(\mathrm{i} \nu m) \tilde{v}_{m, n} e_{m, n}(\mathbf{x}_{h}), \quad \mathcal{D}_{x}^{2} v_{h}:=\sum_{m, n=-M / 2}^{M / 2-1}(\mathrm{i} \nu m)^{2} \tilde{v}_{m, n} e_{m, n}(\mathbf{x}_{h}).
\]
The differential operators $\mathcal{D}_{y}$ and $\mathcal{D}_{y}^{2}$ can be defined in the similar fashion. We define the discrete gradient and Laplacian in the point-wise sense by $\nabla_hv_h := \left(\mathcal{D}_xv_h,\mathcal{D}_yv_h\right)^T$ and $\Delta_hv_h :=\bra{\mathcal{D}_x^2+\mathcal{D}_y^2}v_h,$ respectively.

Denote the space of $L$-periodic grid functions $\mathbb{V}_h:=\{v\mid v=(v_h)$ is $L$-periodic for $\mathbf{x}_h\in\bar{\Omega}_h\}$. For any grid functions $v,w\in\mathbb{V}_h$, we define the discrete inner product $\myinnert{v,w}:=h^2\sum_{\mathbf{x}_h\in{\Omega}_h}v_hw_h$, with the associated $L^2$ norm $\|v\|:=\sqrt{\myinnert{v,v}}$. Also, we will use the maximum norm $\|v\|_{\infty}=\max _{\mathbf{x}_{h} \in \Omega_{h}}|v_{h}|$ and
\[
  \|\nabla_{h} v\|:=\sqrt{h^{2} \sum_{\mathbf{x}_{h} \in \Omega_{h}}|\nabla_{h} v_{h}|^{2}} \quad\text{and}\quad \|\Delta_{h} v\|:=\sqrt{h^{2} \sum_{\mathbf{x}_{h} \in \Omega_{h}}|\Delta_{h} v_{h}|^{2}}.
\]
It is easy to check the discrete Green's formulas, $\langle-\Delta_{h} v, w\rangle=\langle\nabla_{h} v, \nabla_{h} w\rangle$ and $\langle\Delta_{h}^{2} v, w\rangle=\langle\Delta_{h} v, \Delta_{h} w\rangle$, see \cite{GottliebWang:2012,ChengWangWise:2019CCP,ChengWangWiseYue:2016JSC,ShenTangWang:2011book} for more details. Also one has the embedding inequality simulating the Sobolev embedding $H^{2}(\Omega) \hookrightarrow L^{\infty}(\Omega)$,
\begin{align}\label{eq: embedding inequality}
	\|v\|_{\infty} \leq \ck_{\Omega} \brab{\|v\|+\|(I+\Delta_{h}) v\|} \quad\text {for any } v \in \mathbb{V}_{h},
\end{align}
where $\ck_{\Omega}>0$ is only dependent on the domain $\Omega_{h}$. Here and hereafter, any subscripted $\ck$, such as $\ck_{\Omega}, \ck_0, \ck_{\kappa}, \ck_f, \ck_g$ and so on, denotes a fixed constant. The appeared constants may be dependent on the given data  and the solution $\Phi$ but are always independent of the spatial length $h$ and the time-step size $\tau$.

For the underlying volume-conservative PFC model (with $M_h=\Delta_h$), it is also to define a mean-zero function space $\mathring{\mathbb{V}}:=\big\{v\in\mathbb{V} \mid \myinnert{v, 1}=0\big\} \subset \mathbb{V}$. As usual, following the arguments in \cite{ChengWangWise:2019CCP,ChengWangWiseYue:2016JSC}, one can introduce a discrete version of inverse Laplacian operator $(-\Delta_h)^{-\gamma}$ \lan{($\gamma>0$)} as follows. For a grid function $v\in\mathring{\mathbb{V}}$, define
\[
  (-\Delta_h)^{-\gamma} v_h := \sum_{\ell, m=-M/2, (\ell,m) \not=\mathbf{0}}^{M/2-1} \brab{\nu^2(\ell^2+m^2)}^{-\gamma} \tilde{v}_{\ell,m} e_{\ell,m}(\mathbf{x}_h),
\]
and an $H^{-1}$ inner product $\myinnerb{v, w}_{-1} := \myinnerb{(-\Delta_h)^{-1}v, w}$. The associated $H^{-1}$ norm $\mynormb{\cdot}_{-1}$ can be defined by $\mynormb{v}_{-1}:=\sqrt{\myinnerb{v, v}_{-1}}$. 
It is easy to check the following generalized Cauchy-Schwarz inequality by following the proof of \cite[Lemma 2.1]{ChengWangWiseYue:2016JSC},
\begin{align}\label{eq: generalized Cauchy-Schwarz}
	\myinnerb{v, w} \leq \mynormb{\nabla_h v}\mynormb{w}_{-1} \quad\text{for any $v, w\in\mathring{\mathbb{V}}$.}
\end{align}

\subsection{Preliminary lemmas and an equivalent form}

In deriving the discrete energy law \eqref{theoremResult: discrete energy law} from the unified differential form \eqref{eq: unified differential form}, we need the following lemma, which suggests that the uniform boundedness of stage solutions is always necessary.
\begin{lemma}\label{lemma: kappa}
	Assume that the nonlinear bulk $f$ is  continuously differentiable. If  the stabilized parameter $\kappa\geq\max_{\xi_h\in \mathcal{B}_h^n}\mynormt{f^\prime(\xi)}_{\infty}$, then
	\begin{align*}
		\myinnerb{v^n - v^{n-1}, f_{\kappa}(v^{n-1}) - \tfrac12L_{\kappa}(v^n + v^{n-1})} \leq E[v^{n-1}] - E[v^{n}],
	\end{align*}
	 where the function space $\mathcal{B}_h^n:=\big\{\xi_h\in \mathbb{V}_h\,\big|\big.\,\mynormt{\xi}_{\infty}\le \min\{\mynormt{v^{n-1}}_{\infty},\mynormt{v^{n}}_{\infty}\}\big\}.$
\end{lemma}
\begin{proof}
	From the Taylor expansion, we know that
	\begin{align*}
		f(v_h^{n-1})(v_h^n-v_h^{n-1})\leq F(v_h^{n-1}) - F(v_h^n) + \tfrac\kappa2 (v_h^n-v_h^{n-1})^2.
	\end{align*}
	It follows that
	\begin{align*}
		\myinnerb{v^n-v^{n-1} , f_{\kappa}(v^{n-1}) - \tfrac{\kappa}{2}(v^n+v^{n-1})} =&\, \myinnerb{v^n-v^{n-1}, f(v^{n-1})} - \tfrac{\kappa}{2}\|v^n-v^{n-1}\|^2 \\
		\leq&\, \myinnerb{F(v^{n-1}), 1}-\myinnerb{F(v^n), 1}.
	\end{align*}
	Also, due to the symmetry of $L_{h}$, it is easy to get that
	\begin{align*}
		\myinnerb{v^n-v^{n-1}, \tfrac{\kappa}{2}(v^n+v^{n-1}) - \tfrac{1}{2}L_{\kappa}(v^n+v^{n-1})}=&\, - \tfrac{1}{2} \myinnerb{v^n-v^{n-1}, L_{h}(v^n+v^{n-1})} \\
		 = &\,\tfrac{1}{2}\myinnerb{v^{n-1}, L_{h}v^{n-1}} - \tfrac{1}{2}\myinnerb{v^n, L_{h}v^n}.
	\end{align*}
	From the definition \eqref{eq: energy functional}, one has $E[v^{n}]=\tfrac{1}{2}\myinnerb{v^n, L_{h}v^n}+\myinnerb{F(v^n), 1}$. Thus adding up the above two results yields the claimed inequality and completes the proof.
\end{proof}

As done in \cite{LiQiaoWangZheng:2024arxiv}, we will take the advantage of the inherent $H^2$ norm boundedness in the associated free energy functional \eqref{eq: energy functional}, that is, the boundedness of original energy functional $E[u^{n,i}]$ implies the $H^2$ boundedness of stage solutions $u_h^{n,i}$, as stated in the next lemma, in which some norms of the nonlinear term $f$ are also bounded. 
\begin{lemma}\label{lemma: Energy implies H2 norm}
	If $E[v] \leq E[\phi^0]$ for any $v\in\mathbb{V}_h$, then
	\[
	\mynormb{(I+\Delta_h) v} + \mynormb{v} \leq \ck_0=\sqrt{4E[\phi^{0}]+\mynormt{1+\varepsilon}^2}.
	\]
	Furthermore, for the nonlinear bulk $f(v)=\varepsilon v-v^3$, it holds that
	\begin{align*}
	\max\big\{\mynormb{f(v)}, \mynormb{f_{\kappa}(v)}\big\} \leq \ck_{\kappa} \mynormb{v}\;\;\text{and}\;\;
	\mynormb{\nabla_h f(v)} \leq \ck_{g}\bra{\mynormb{(I+\Delta_h) v} + \mynormb{v}}.
	\end{align*}
\end{lemma}
\begin{proof}
  The fact $(a^2-1-\varepsilon)^2\ge0$ gives $\myinnerb{v^4-2\varepsilon v^2,1}\ge 2\mynormt{v}^2-\mynormt{1+\varepsilon}^2$. Recalling the definition $E[v]$, one has $$\ck_0^2=4E[\phi^0]+\mynormt{1+\varepsilon}^2\ge 2\mynormt{(I+\Delta_h)v}^2+2\mynormt{v}^2, $$ which arrives at the first result, also see \cite{SunZhaoetal:2022CNSNS}. It implies that $\mynormt{v}_{\infty}\le \ck_{\Omega}\ck_0$ due to the embedding inequality \eqref{eq: embedding inequality}. 
	Thus we have $\mynormb{v^3}\le \ck_{\Omega}^2\ck_0^2\mynormb{v}$
	such that $$\max\big\{\mynormb{f(v)}, \mynormb{f_{\kappa}(v)}\big\}\le \ck_{\kappa}\mynormb{v}$$  by taking the constant $\ck_{\kappa}:=\varepsilon+\kappa+\ck_{\Omega}^2\ck_0^2$. 
	By using \cite[Lemma 1]{GottliebWang:2012}, one has
	$$\mynormb{\nabla_hv^3}=\mynormb{\nabla(I_Mv^3)}_{L^2}\le 
	2\mynormb{\nabla v^3}_{L^2}=6\mynormb{v^2\nabla v}_{L^2}.$$ Thus, with the help of the Sobolev embedding inequality, there exists a positive constant $\ck_H$ such that $\mynorm{\nabla_hv^3}\le \ck_H\ck_0^2\mynorm{\nabla_h v}$. Thus the last result follows from the  fact $\mynorm{\nabla_h v}\le \mynorm{(I+\Delta_h) v} + \mynorm{v}$ with $\ck_g:=\ck_H\ck_0^2$.
  The proof is completed.
\end{proof}

As seen, the key point to our aim is to establish the maximum norm boundedness of the numerical solution at all stages, which requires the $L^2$ norm and $H^2$ semi-norm bounds  according to  Lemma \ref{lemma: Energy implies H2 norm}. In the following, we will derive the $H^2$ semi-norm boundedness of the stage solution via the differential form \eqref{eq: unified differential form}; while the $L^2$ norm boundedness will be evaluated via its equivalent DOC form, which  is formulated by using the so-called discrete orthogonal convolution (DOC) kernels \cite{LiLiao:2022SINUM,LiTangZhou:2024SCM,LiaoZhang:2021}.
\lan{For the kernels $d_{ij}(z)$, the DOC kenrels $\theta_{kj}(z)$ is defined by, for each index $k=1,2,\cdots, s$,}
\lan{\begin{equation*}
  \theta_{kk}(z) := d_{kk}^{-1}(z) \quad\text{and}\quad \theta_{kj}(z) := -d_{jj}^{-1}(z)\sum_{\ell=j+1}^k \theta_{k\ell}(z)d_{\ell j}(z) \quad\text{for $1 \leq j \leq k-1$.}
\end{equation*}
It is easy to check the following discrete orthogonal identity, 
\begin{equation*}
  \sum_{\ell=j}^m \theta_{m\ell}(z) d_{\ell j}(z) \equiv \delta_{m j} \quad\text{for $1\leq j \leq m \leq s$,}
\end{equation*}
where $\delta_{m j}$ is the Kronecker delta symbol with $\delta_{m j} = 0$ if $j\not=m$.} 
Multiplying the differential form \eqref{eq: unified differential form} by the DOC kernels $\theta_{ki}$ and summing $i$ from 1 to $k$, one can exchange the summation order and apply the above orthogonal identity to obtain the following DOC form
\begin{equation}\label{eq: DOC form}
	\delta_{\tau} u_h^{n, i+1} = \tau \sum_{j=1}^i \theta_{ij}(\tau M_h L_{\kappa}) M_h \kbraB{L_{\kappa} u_h^{n, j+\frac12} - f_{\kappa}(u_h^{n, j})}
\end{equation}
for $n\ge1$ and $1 \leq i \leq s$.
Interested readers can find the related researches \cite{LiaoZhang:2021,LiaoJiZhang:2022IMA,LiaoKang:2023IMA,LiLiao:2022SINUM} on the $L^2$ norm analysis of time-stepping methods using the backward differentiation formulas, where the differential form of a numerical scheme is always \lan{reformulated} into its equivalent DOC form, and the $L^2$ norm estimate of the numerical solution can be performed via  the standard discrete energy arguments. 

Actually, the above discrete orthogonal identity says that the DOC matrix $$\Theta(z)=\kbrab{\theta_{ij}(z)}_{i,j=1}^s=D(z)^{-1}.$$ Thus the DOC matrix $\Theta(z)$ is positive \lan{semi-}definite if the differentiation matrix $D(z)$ is positive definite. It arrives at the following result.
\begin{lemma}\label{lemma: bound quadratic form}	
Assume that the differentiation matrix $D(z)$ is positive definite. If there exists a finite constant $\lambda_{\es}>0$, independent of $z$, such that it is not larger than any minimum eigenvalues of $\mathcal{S}[D(z)]$ for $z\le0$, then for any time sequences $\{v^{i},u^i \mid i \geq 1\}$, it \lan{holds} that
	\begin{align*}
		&\,\sum_{i=1}^{k} \sum_{j=1}^{i} \myinnerb{d_{ij}(z) v^{j}, v^{i}} \geq \lambda_{\es}\sum_{i=1}^{k}\mynormb{v^{i}}^2\quad\text{for $1\le k\le s$,}\\
	&\,\sum_{i=1}^{k} \sum_{j=1}^{i} \myinnerb{\theta_{ij}(z) v^{j}, u^{i}} \leq \frac{1}{\lambda_{\es}} \sum_{i=1}^{k}\mynormb{v^{i}}\mynormb{u^{i}}\quad\text{for $1\le k\le s$.}
	\end{align*}
\end{lemma}
\begin{proof}The first inequality is obvious since the matrix $D-\lambda_{\es}I$ is positive semi-definite, while the second result follows from an upper bound estimate of the spectral norm of the DOC matrix $\Theta$. Actually,  
the maximum eigenvalue of $\Theta^T\Theta$ can be bounded by $1/\lambda_{\es}$, see \cite[Lemma 5.2]{Chenetal:2024JCM}. It completes the proof.
	\end{proof}

We are ready to present the complete proof of Theorem \ref{theorem: energy stability}. To make our proof more clear, we split it into two parts: (Part A) \cref{subsection: SH} addresses the case $M_h=-I$ corresponding to the SH equation; (Part B) \cref{subsection: PFC} discusses the case $M_h=\Delta_h$ corresponding to the PFC equation. 

\subsection{Proof of Theorem \ref{theorem: energy stability} for  $M_h=-I$}
\label{subsection: SH}

\begin{proof}(\textit{Part A: $M_h=-I$}) If $\mynormt{u^{n,i}}_\infty$ is bounded with a proper bound $\tau_0$ of the time-step size, it is easy to determine the constant $\ck_f$ via Lemma \ref{lemma: kappa} and establish the original energy law \eqref{theoremResult: discrete energy law} via a standard energy argument, that is, make the inner product of \eqref{eq: unified differential form} with $2M_h^{-1}\delta_{\tau}u_h^{n, i+1}$ and sum $i$ from $i=1$ to $k$, cf. \cite[Theorem 2.1]{LiaoWang:2024MCOM}.
	
The main task of the proof is to verify the uniform boundedness of solutions $u_h^{n,i}$ at all stages by using the mathematical induction for the original energy bound
\begin{equation}\label{SH: induction aim}
  E[u^{n, i}] \leq E[\phi^0]\quad\text{for $n\ge1$ and $1\le i\le s$.}
\end{equation}
	
Obviously, the inequality \eqref{SH: induction aim} holds for $n=1$ and $i=1$ since $E[u^{1,1}] = E[\phi^0]$.
	We put the inductive hypothesis 
	\begin{equation*}
	E[u^{\ell, j}] \leq E[\phi^0]\quad\text{for $1 \leq \ell \leq n$ and $1 \leq j \leq k$,}
	\end{equation*}
	that is, the energy bound \eqref{SH: induction aim} holds for $1 \leq \ell \leq n$ and $1 \leq j \leq k$.
	Lemma \ref{lemma: Energy implies H2 norm} gives
	\begin{equation}\label{SH: assumption}
		\mynormb{(I+\Delta_h)u^{\ell,j}} + \mynormb{u^{\ell,j}} \leq \ck_0 \quad \text{and} \quad \mynormb{u^{\ell,j}}_\infty \leq \ck_{\Omega}\ck_0
	\end{equation}
	for $1 \leq \ell \leq n$ and $1 \leq j \leq k$, according to the embedding inequality \eqref{eq: embedding inequality}. It is to prove that $E[u^{n, k+1}] \leq E[\phi^0]$ by the $L^2$ norm and $H^2$ semi-norm estimates.
	
	(\textit{$L^2$ norm estimate})
	Making the inner product of the DOC form \eqref{eq: DOC form} with $2u^{n, i+\frac12}$, applying the discrete \lan{Green's} formulas and summing $i$ from 1 to $k$ yield that
	\begin{align}\label{SH pf: L2 inner product}
		\mynormb{u^{n,k+1}}^2 - \mynormb{u^{n,1}}^2
	=2\tau \sum_{i=1}^k \sum_{j=1}^i \myinnerB{\theta_{ij} \kbrab{f_{\kappa}(u^{n, j})-L_{\kappa} u^{n, j+\frac12}}, u^{n, i+\frac12}}, 
	\end{align}
	where we used the fact that $(a+b)(a-b) = a^2-b^2$. Because the DOC matrix $\Theta(z)$ is positive \lan{semi-}definite, the linear term at the right side of \eqref{SH pf: L2 inner product} reads
	\[
	-2\tau \sum_{i=1}^k \sum_{j=1}^i \myinnerB{\theta_{ij}(\tau M_hL_{\kappa}) L_{\kappa} u^{n, j+\frac12}, u^{n, i+\frac12}} \leq 0.
	\]
	With the help of the hypothesis \eqref{SH: assumption} and Lemma \ref{lemma: Energy implies H2 norm}, the nonlinear term at the right side of \eqref{SH pf: L2 inner product} can be bounded by Lemma \ref{lemma: bound quadratic form} as follows
	\begin{align*}
		2\tau \sum_{i=1}^k \sum_{j=1}^i \myinnerB{\theta_{ij}(\tau M_hL_{\kappa}) f_{\kappa}(u^{n, j}), u^{n, i+\frac12}} \leq \frac{2\ck_{\kappa} \tau}{\lambda_{\es}} \sum_{i=1}^k \mynormb{u^{n, i}}\mynormb{u^{n, i+\frac12}}.
	\end{align*}
	Substituting the above two inequalities into \eqref{SH pf: L2 inner product}, one has
	\begin{align*}
		\mynormb{u^{n,k+1}}^2\le \mynormb{u^{n,1}}^2 + \frac{2\ck_{\kappa}\tau}{\lambda_{\es}} \sum_{i=1}^k \mynormb{u^{n, i}} \mynormb{u^{n, i+\frac12}}.
	\end{align*}
	For any $k$, we can find some $k_{0}$ $(0\le k_0\le k)$ such that $\|u^{n,k_{0}+1}\| = \max_{0\leq i \leq k}\|u^{n,i+1}\|$. By taking $k=k_0$ in the above inequality, one can obtain
	\begin{align*}
		\mynormb{u^{n,k_0+1}}^2 \leq \mynormb{u^{n,1}} \mynormb{u^{n,k_0+1}} + \frac{2\ck_{\kappa}\tau}{\lambda_{\es}} \sum_{i=1}^{k_0} \mynormb{u^{n,i}} \mynormb{u^{n,k_0+1}},
	\end{align*}
	that is, by using the hypothesis \eqref{SH: assumption},
	\begin{equation*}
		\mynormb{u^{n,k+1}} \leq \|u^{n,k_{0}+1}\| \leq \mynormb{u^{n,1}} + \frac{2\ck_{\kappa} \tau}{\lambda_{\es}} \sum_{i=1}^k \mynormb{u^{n,i}} \leq \ck_0 + \frac{2 \ck_{\kappa} s \tau}{\lambda_{\es}} \ck_0.
	\end{equation*}
	It gives $\mynormt{u^{n,k+1}} \leq 2\ck_0$ if the time-step size  $\tau\leq \lambda_{\es} / (2\ck_{\kappa} s)$.
	
	(\textit{$H^2$ semi-norm estimate}) Making the inner product of the differential form \eqref{eq: unified differential form} with $2M_h^{-1}\delta_{\tau}u^{n, i+\frac12}$, applying the discrete \lan{Green's} formulas  and summing the stage index $i$ from 1 to $k$, one can derive that 
	\begin{align}\label{SH pf: H2 inner product}
		&\,\mynormb{(I+\Delta_h)u^{n,k+1}}^2 - \mynormb{(I+\Delta_h)u^{n,1}}^2 +\kappa\sum_{i=1}^k \mynormb{\delta_{\tau}u^{n, i+1}}^2 \\
		=&\, -\frac{2}{\tau} \sum_{i=1}^k \sum_{j=1}^i \myinnerB{d_{ij}\delta_{\tau} u^{n, j+1}, \delta_{\tau}u^{n, i+1}} + 2\sum_{i=1}^k \myinnerb{f(u^{n, i}), \delta_{\tau}u^{n, i+1}}, \notag
	\end{align}
	where the definitions of $L_{\kappa}$ and $f_{\kappa}$ have been used.
	The positive definiteness of the differentiation matrix $D(z)$, see Lemma \ref{lemma: bound quadratic form}, yields
	\[
  	-\frac2\tau \sum_{i=1}^k \sum_{j=1}^i \myinnerB{d_{ij}(\tau M_hL_{\kappa}) \delta_{\tau} u^{n, j+1}, \delta_{\tau}u^{n, i+1}} \leq -\frac{2\lambda_{\es}}\tau \sum_{i=1}^k \mynormb{\delta_{\tau}u^{n, i+1}}^2.
	\]
	With the help of the hypothesis \eqref{SH: assumption} and Lemma \ref{lemma: Energy implies H2 norm}, the nonlinear term at the right side of \eqref{SH pf: H2 inner product} can be bounded by 
		\begin{align*}
		2\sum_{i=1}^k \myinnerb{f(u^{n, i}), \delta_{\tau}u^{n, i+1}}
		\le&\, \frac{2\lambda_{\es}}\tau \sum_{i=1}^k \mynormb{\delta_{\tau}u^{n, i+1}}^2
		+\frac{\ck_{\kappa}^2\tau}{2\lambda_{\es}} \sum_{i=1}^k \mynormb{u^{n, i}}^2.
	\end{align*}
	Reminding the hypothesis bounds in \eqref{SH: assumption}, we substitute the above two estimates into \eqref{SH pf: H2 inner product} and find that
	\begin{align*}
		\mynormb{(I+\Delta_h)u^{n,k+1}}^2 \leq&\, \mynormb{(I+\Delta_h)u^{n,1}}^2 + \frac{\ck_{\kappa}^2\tau}{2\lambda_{\es}} \sum_{i=1}^k \mynormb{u^{n,i}}^2 \leq \ck_0^2 + \frac{\ck_{\kappa}^2s\tau}{2\lambda_{\es}} \ck_0^2.
	\end{align*}
	One has the $H^2$ semi-norm estimate $\mynormt{(I+\Delta_h)u^{n,k+1}} \leq 2\ck_0$ if
$\tau\leq 6\lambda_{\es} / (\ck_{\kappa}^2s)$.
	
	(\textit{Maximum norm estimate}) Adding up the above $L^2$ norm and $H^2$ semi-norm estimates gives ($\ck_0$ is independent of the stabilized parameter $\kappa$)
	\begin{align*}
		\mynormt{u^{n,k+1}} + \mynormt{(I+\Delta_h)u^{n,k+1}} \leq 4\ck_0
		\quad\text{if $\tau\le \min\{\lambda_{\es} / (2\ck_{\kappa} s),6\lambda_{\es} / (\ck_{\kappa}^2s)\}$. }
	\end{align*}
	By taking a constant $\tau_0:=\min\{\lambda_{\es} / (2\ck_{\kappa} s),6\lambda_{\es} / (\ck_{\kappa}^2s)\}$, which may be dependent on the parameter $\kappa$,  we arrive at the maximum norm bound, $\mynormt{u^{n,k+1}}_\infty\leq4\ck_{\Omega} \ck_0$, according to the embedding inequality \eqref{eq: embedding inequality}. Thus by taking the constant $$\ck_{f}:=\max_{\mynormt{\xi}_{\infty}\le 4\ck_{\Omega} \ck_0}\mynormt{f^\prime(\xi)}_{\infty},$$ one can apply Lemma \ref{lemma: kappa} to recover the original energy law \eqref{theoremResult: discrete energy law} at the stage $t_{n,k+1}$, which gives the desired energy bound $E[u^{n, k+1}] \leq E[u^{n, 1}] \leq E[\phi^0]$
	and completes the mathematical induction. 
	
	Therefore, the original energy bound \eqref{SH: induction aim} holds for any $n\ge1$ and $1\le i\le s$. As a by-product result, one can apply Lemma \ref{lemma: Energy implies H2 norm} to obtain the following global-in-time estimates of stage solutions, $\mynormb{(I+\Delta_h)u^{n,i}} + \mynormb{u^{n,i}}  \leq \ck_0$ for $n\ge1$ and $1\le i\le s$.
	It completes the proof.
\end{proof}

\subsection{Proof of Theorem \ref{theorem: energy stability} for $M_h=\Delta_h$}
\label{subsection: PFC}

\begin{proof}
  (\textit{Part B: $M_h=\Delta_h$}) In this case, it is easy to check that the stage solutions $u^{n,i}$ preserve the volume, $\myinnerb{u^{n,i}, 1}=\myinnerb{u^{1,1}, 1}$ for $n\ge1$ and $1\le i\le s$.
As similar to the case  (\textit{Part A: $M_h=-I$}),  the main task is to verify the uniform boundedness of stage solutions $u_h^{n,i}$ by using the mathematical induction for the  energy bound
  \begin{equation}\label{PFC: induction aim}
    E[u^{n, i}] \leq E[\phi^0]\quad\text{for $n\ge1$ and $1\le i\le s$.}
  \end{equation}

  Obviously, the energy bound \eqref{PFC: induction aim} holds for $n=1$ and $i=1$ since $E[u^{1,1}] = E[\phi^0]$.
  We propose the following inductive hypothesis
  \begin{equation*}
    E[u^{\ell, j}] \leq E[\phi^0]\quad\text{for $1 \leq \ell \leq n$ and $1 \leq j \leq k$,}
  \end{equation*}
  that is, the energy bound \eqref{PFC: induction aim} holds for $1 \leq \ell \leq n$ and $1 \leq j \leq k$.
  Lemma \ref{lemma: Energy implies H2 norm} gives
  \begin{equation}\label{PFC: assumption}
    \mynormb{(I+\Delta_h)u^{\ell,j}} + \mynormb{u^{\ell,j}}  \leq \ck_0 \quad \text{and} \quad \mynormb{u^{\ell,j}}_\infty \leq \ck_{\Omega}\ck_0
  \end{equation}
  for $1 \leq \ell \leq n$ and $1 \leq j \leq k$, according to the embedding inequality \eqref{eq: embedding inequality}. It is to prove that $E[u^{n, k+1}] \leq E[\phi^0]$ by the $H^2$ semi-norm and $L^2$ norm estimates.

  (\textit{$H^2$ semi-norm estimate}) Making the inner product of the differential form \eqref{eq: unified differential form} with $2M_h^{-1}\delta_{\tau}u^{n, i+\frac12}$, applying the discrete \lan{Green's} formulas and summing the stage index $i$ from 1 to $k$, one can derive that
    \begin{align}\label{PFC pf: H2 inner product}
    &\mynormb{(I+\Delta_h)u^{n,k+1}}^2 - \mynormb{(I+\Delta_h)u^{n,1}}^2 +\kappa\sum_{i=1}^k \mynormb{\delta_{\tau}u^{n, i+1}}^2 \\
    &=\frac{2}{\tau} \sum_{i=1}^k \sum_{j=1}^i \myinnerB{d_{ij} \delta_{\tau} u^{n, j+1}, M_h^{-1}\delta_{\tau}u^{n, i+1}}+ 2\sum_{i=1}^k \myinnerb{f(u^{n, i}), \delta_{\tau}u^{n, i+1}}. \notag
  \end{align}
  The positive semi-definiteness of the matrix $D-\lambda_{\es}I$, see Lemma \ref{lemma: bound quadratic form}, yields
  \[
    \frac2\tau \sum_{i=1}^k \sum_{j=1}^i \myinnerB{d_{ij}(\tau M_h L_{\kappa}) \delta_{\tau} u^{n, j+1}, M_h^{-1}\delta_{\tau}u^{n, i+1}} \leq -\frac{2\lambda_{\es}}\tau \sum_{i=1}^k \mynormb{\delta_{\tau}u^{n, i+1}}_{-1}^2.
  \]
  With the help of the generalized Cauchy-Schwarz inequality \eqref{eq: generalized Cauchy-Schwarz}, the nonlinear term at the right side of \eqref{PFC pf: H2 inner product} can be bounded by 
  \begin{align*}
     &\,2\sum_{i=1}^k \myinnerb{f(u^{n, i}), \delta_{\tau}u^{n, i+1}}
   \le 2 \sum_{i=1}^k \mynormb{\delta_{\tau}u^{n, i+1}}_{-1} \mynormb{\nabla_h f(u^{n, i})}\\
    &\,\le \frac{2\lambda_{\es}}\tau \sum_{i=1}^k \mynormb{\delta_{\tau}u^{n, i+1}}_{-1}^2
    +\frac{\ck_{g}^2\tau}{2\lambda_{\es}} \sum_{i=1}^k\brab{\mynormb{(I+\Delta_h)u^{n,i}} + \mynormb{u^{n,i}}}^2\\
    &\,\le \frac{2\lambda_{\es}}\tau \sum_{i=1}^k \mynormb{\delta_{\tau}u^{n, i+1}}_{-1}^2
    +\frac{\ck_{g}^2s\tau}{2\lambda_{\es}} \ck_0^2,
  \end{align*}
  where Lemma \ref{lemma: Energy implies H2 norm} together with the inductive hypothesis \eqref{PFC: assumption} was applied in the last inequality.
  Thus we substitute the above two estimates into \eqref{PFC pf: H2 inner product} and find that
  \begin{align*}
    \mynormb{(I+\Delta_h)u^{n,k+1}}^2 \leq \ck_0^2 + \frac{\ck_{g}^2s\tau}{2\lambda_{\es}} \ck_0^2.
  \end{align*}
  One has the $H^2$ semi-norm estimate, $\mynormb{(I+\Delta_h)u^{n,k+1}} \leq 2\ck_0$ if $\tau\leq 6\lambda_{\es} / (\ck_{g}^2s)$.

  (\textit{$L^2$ norm estimate})
  Making the inner product of the DOC form \eqref{eq: DOC form} with $2u^{n, i+\frac12}$, applying the discrete \lan{Green's} formulas and summing the stage index $i$ from 1 to $k$, one can obtain that
  \begin{align}\label{PFC pf: L2 inner product}
    \mynormb{u^{n,k+1}}^2 - \mynormb{u^{n,1}}^2
    =2\tau \sum_{i=1}^k \sum_{j=1}^i \myinnerB{\theta_{ij}  \kbrab{L_{\kappa} u^{n, j+\frac12}-f_{\kappa}(u^{n, j})}, \Delta_hu^{n, i+\frac12}}. 
  \end{align}
  Because the DOC matrix $\Theta(z)$ is positive \lan{semi-}definite, the linear term at the right side of the equality \eqref{PFC pf: L2 inner product} can be dropped, that is,
  \[
    2\tau \sum_{i=1}^k \sum_{j=1}^i \myinnerB{\theta_{ij}(\tau M_h L_{\kappa})  L_{\kappa} u^{n, j+\frac12}, \Delta_hu^{n, i+\frac12}} \leq 0.
  \]
  With the help of Lemma \ref{lemma: Energy implies H2 norm} and the hypothesis \eqref{PFC: assumption}, the nonlinear term at the right side of \eqref{PFC pf: L2 inner product} can be bounded by Lemma \ref{lemma: bound quadratic form} as follows
  \begin{align*}
    &\, -2\tau \sum_{i=1}^k \sum_{j=1}^i \myinnerB{\theta_{ij}(\tau M_h L_{\kappa}) f_{\kappa}(u^{n, j}), \Delta_hu^{n, i+\frac12}} \\
    =&\, -2\tau \sum_{i=1}^k \sum_{j=1}^i \myinnerB{\theta_{ij}(\tau M_h L_{\kappa}) f_{\kappa}(u^{n, j}), (I+\Delta_h) u^{n, i+\frac12}-u^{n, i+\frac12}} \\
    \leq&\, \frac{2\ck_{\kappa} \tau}{\lambda_{\es}} \sum_{i=1}^k \mynormb{u^{n, i}} \brab{\mynormb{(I+\Delta_h) u^{n, i+\frac12}} + \mynormb{u^{n, i+\frac12}}}.
  \end{align*}
  Substituting the above two inequalities into \eqref{PFC pf: L2 inner product}, one has
  \begin{align*}
    \mynormb{u^{n,k+1}}^2\le \mynormb{u^{n,1}}^2 + \frac{2\ck_{\kappa} \tau}{\lambda_{\es}} \sum_{i=1}^k \mynormb{u^{n, i}} \brab{\mynormb{(I+\Delta_h) u^{n, i+\frac12}} + \mynormb{u^{n, i+\frac12}}}.
  \end{align*}
  For any $k$, we can find some $k_{0}$ $(0\le k_0\le k)$ such that $\|u^{n,k_{0}+1}\| = \max_{0\leq i \leq k}\|u^{n,i+1}\|$. By taking $k=k_0$ in the above inequality, one can obtain
  \begin{align*}
    \mynormb{u^{n,k_0+1}}^2 \leq \mynormb{u^{n,1}} \mynormb{u^{n,k_0+1}} + \frac{2\ck_{\kappa} \tau}{\lambda_{\es}} \sum_{i=1}^k \mynormb{u^{n,k_0+1}} \brab{\mynormb{(I+\Delta_h) u^{n, i+\frac12}} + \mynormb{u^{n, i}}}.
  \end{align*}
  Then, by using the triangle inequality and the hypothesis in \eqref{PFC: assumption},
  \begin{equation*}
    \mynormb{u^{n,k+1}} \leq \mynormb{u^{n,1}} + \frac{2\ck_{\kappa} \tau}{\lambda_{\es}} \sum_{i=1}^k \brab{\mynormb{(I+\Delta_h) u^{n, i+\frac12}} + \mynormb{u^{n, i}}} \leq \ck_0 + \frac{4 \ck_{\kappa} s \tau}{\lambda_{\es}} \ck_0,
  \end{equation*}
  where the available $H^2$ semi-norm estimate, $\mynormb{(I+\Delta_h)u^{n,k+1}} \leq 2\ck_0$, was also used. It arrives at the  $L^2$ norm estimate, $\mynormb{u^{n,k+1}} \leq 2\ck_0$ if $\tau\leq \lambda_{\es} / (4\ck_{\kappa} s)$.
	
	(\textit{Maximum norm estimate}) Adding up the above $L^2$ norm and $H^2$ semi-norm estimates gives ($\ck_0$ is independent of the stabilized parameter $\kappa$)
  \begin{align*}
    \mynormt{u^{n,k+1}} + \mynormt{(I+\Delta_h)u^{n,k+1}} \leq 4\ck_0
    \quad\text{if $\tau\le \min\{\lambda_{\es} / (4\ck_{\kappa} s),6\lambda_{\es} / (\ck_{g}^2s)\}$. }
  \end{align*}
  By taking a constant $\tau_0:=\min\{\lambda_{\es} / (4\ck_{\kappa} s),6\lambda_{\es} / (\ck_{g}^2s)\}$, which may be dependent on the parameter $\kappa$,  we arrive at the maximum norm bound, $\mynormt{u^{n,k+1}}_\infty\leq4\ck_{\Omega} \ck_0$, according to \eqref{eq: embedding inequality}. Thus by taking the constant $\ck_{f}:=\max_{\mynormt{\xi}_{\infty}\le 4\ck_{\Omega} \ck_0}\mynormt{f^\prime(\xi)}_{\infty},$ one can apply Lemma \ref{lemma: kappa} to recover the original energy law \eqref{theoremResult: discrete energy law} at the stage $t_{n,k+1}$ by making the inner product of the differential form \eqref{eq: unified differential form} with $2M_h^{-1}\delta_{\tau}u_h^{n, i+1}$ and summing $i$ from $i=1$ to $k$. It gives the desired bound $E[u^{n, k+1}] \leq E[u^{n, 1}] \leq E[\phi^0]$
  and completes the mathematical induction. 
  
  Therefore, the original energy bound \eqref{PFC: induction aim} holds for any $n\ge1$ and $1\le i\le s$. As a by-product result, one can apply Lemma \ref{lemma: Energy implies H2 norm} to obtain the following global-in-time estimates of stage solutions,
 $\mynormb{(I+\Delta_h)u^{n,i}} + \mynormb{u^{n,i}}  \leq \ck_0$ for $n\ge1$ and $1\le i\le s$.
 It completes the proof.
\end{proof}

To end this section, we remark that, Theorem \ref{theorem: energy stability} makes the associated convergence analysis to be somewhat easy because the nonlinear error control equation can be reduced into a linear one essentially and we omit it here. 
Interesting readers can find the related discussions in the literature including \cite{LiQiao:2024SISC,SunZhangetal:2024arxiv}.
Furthermore, with the help of the unified differential form \eqref{eq: unified differential form} and the standard energy arguments, the present unified theory make the full use of the global property (positive definiteness) of RK coefficients rather than the local properties of RK coefficients so that it is much conciser than the recent technical estimates in \cite{LiQiaoWangZheng:2024arxiv,SunZhangetal:2024arxiv} on the global-in-time energy estimates of specific EERK methods.

It is to mention that, the present stability analysis would be closely related to the so-called internal stability analysis for RK methods, which would be useful to control the stability associated with each stage in addition to each step beyond traditional step-wise stability. The discrete energy arguments of the internal stability can be found in \cite{Hairer:1984book,KennedyCarpenter:2016}, in which a stage-wise analog to the algebraic-stability matrix of algebraically stable implicit RK methods was considered. Utilizing the unified differential form \eqref{eq: unified differential form} and the associated DOC form \eqref{eq: DOC form} of the IERK, EERK and CIFRK methods, our theoretical approach can be regarded as a novel approach for the internal nonlinear stability of some efficient (not necessarily algebraically stable) RK methods  for dissipative semilinear parabolic problems.

\section{Lower bounds $\lambda_{\es}$ for IERK, EERK and CIFRK methods}\label{section: lambda}


From the stability analysis in \cref{section: theoretical framework},  one can see that, the condition of Lemma \ref{lemma: bound quadratic form} is essential to the result of Theorem \ref{theorem: energy stability}. For different type RK methods having positive definite differentiation matrices $D(z)$, the common condition is that there exists a finite constant $\lambda_{\es}>0$, independent of $z$, such that it is not larger than any minimum eigenvalues of $\mathcal{S}[D(z)]$ for $z<0$.  
We  will evaluate the existence of lower bounds $\lambda_{\es}$ by revisiting some IERK, EERK and CIFRK methods in the literature. 
Note that, for several parameterized RK methods, such as IERK methods in \cite{AscherRuuthSpiteri:1997,LiaoWangWen:2024IERK} or EERK methods in \cite{HochbruckOstermann:2005SINUM,HochbruckOstermann:2010EERK,LiaoWang:2024MCOM}, we consider some specific RK methods with the parameter chosen by the average dissipation rate, which was proposed recently in \cite{LiaoWang:2024MCOM,LiaoWangWen:2024IERK} to roughly measure the degree of preservation of original energy dissipation law.



\subsection{IERK methods}

In general, the $(s+1)$-stage IERK methods that are globally stiffly-accurate have the following general form \cite{AscherRuuthSpiteri:1997,CooperSayfy:1983,FuTangYang:2024MCOM,KennedyCarpenter:2003ANM,LiaoWangWen:2024IERK}
\begin{align}\label{IERK: general form}     
  u^{n,i} = u^{n,1} + \tau\sum_{j=1}^i a_{ij} M_h L_\kappa u^{n,j} - \tau\sum_{j=1}^{i-1} \hat{a}_{ij} M_h f_\kappa(u^{n,j})
\end{align}
for $n\ge1$ and $1\le i\le s+1$, with the following Butcher tableaux
\begin{equation*}
  \begin{array}{c|ccccc|ccccc}
    c_{1} & 0 &  &  &  & & 0   \\
    c_{2} & a_{21} & a_{22} &  &  &  & \hat{a}_{21} & 0  \\
    c_{3} & a_{31} & a_{32} & a_{33} &  &  & \hat{a}_{31} & \hat{a}_{32} & 0  \\
    \vdots & \vdots & \vdots & \ddots & \ddots &  \\[2pt]
    c_{s+1} & a_{s+1,1} & a_{s+1,2} &  \cdots  & a_{s+1,s}   & a_{s+1,s+1}  & \hat{a}_{s+1,1} & \hat{a}_{s+1,2} &  \cdots  & \hat{a}_{s+1,s}   & 0 \\[2pt]
    \hline  & a_{s+1,1} & a_{s+1,2} &  \cdots  & a_{s+1,s}   & a_{s+1,s+1} & \hat{a}_{s+1,1} & \hat{a}_{s+1,2} &  \cdots  & \hat{a}_{s+1,s}   & 0
  \end{array},
\end{equation*} 
where the abscissas $c_1=0$ and $c_{s+1}=1$. Here the  coefficients $a_{ij}$ and $\hat{a}_{ij}$ are real numbers that can be dependent on some given \lan{parameters} (see the parameterized IERK methods in \cite{AscherRuuthSpiteri:1997,LiaoWangWen:2024IERK}), but always be independent of any discrete operators including $M_h$ and $L_\kappa$. By the canopy node condition $\sum_{j=1}^i a_{ij} = c_i = \sum_{j=1}^{i-1} \hat{a}_{ij}$, one can reformulate the  form \eqref{IERK: general form} into the steady-state preserving form \cite[Eq. (2.4)]{LiaoWangWen:2024IERK}
\begin{align}\label{IERK: steady-state preserving form}
 u_h^{n, i+1} =u_h^{n,1} &\,+ \tau\sum_{j=1}^{i} a_{i+1,j+1}M_h \kbrab{L_\kappa u_h^{n,j+1}-L_\kappa u_h^{n,1}}\\
&\,- \tau\sum_{j=1}^i \hat{a}_{i+1,j} M_h \kbrab{f_\kappa(u_h^{n,j})-L_\kappa u_h^{n,1}}  
\quad\text{for $n\ge1$ and $1\le i\le s$,}\notag
\end{align}
where the stiff terms with the coefficients $a_{i,1}$ for $2\leq i \leq s+1$ are dropped. In this sense, we define the lower triangular coefficient matrices $A_{\mathrm{I}}:=\brab{a_{i+1,j+1}}_{i,j=1}^{s}$ and $A_{\mathrm{E}}:=\brab{\hat{a}_{i+1,j}}_{i,j=1}^{s}$ for the implicit and explicit parts, respectively. 

Let $E_{s}:=(1_{i\ge j})_{s\times s}$ be the lower triangular matrix full of element 1 and $I_{s}$ be the identity matrix of the same size as $A_{\mathrm{I}}$. With the help of \eqref{IERK: steady-state preserving form}, the IERK method \eqref{IERK: general form}  can be formulated into the differential form \eqref{eq: unified differential form} with the associated differentiation matrix $D=\kbrab{d_{ij}(z)}_{s\times s}$ defined by, cf. \cite[subsection 2.1]{LiaoWangWen:2024IERK},
\begin{align*}\label{IERK: differentiation matrix}     
  D(z):=D_{\mathrm{E}}-zD_{\mathrm{EI}} \quad\text{for $z<0$,}
\end{align*}
where $D_{\mathrm{E}}:=\kbrab{d_{ij}^{(e)}}_{s\times s}=A_{\mathrm{E}}^{-1}E_{s}$ and $D_{\mathrm{EI}}:=A_{\mathrm{E}}^{-1}A_{\mathrm{I}}E_{s}-E_{s}+\tfrac{1}2I_{s}$. 
One can claim that \cref{theorem: energy stability} holds for any IERK methods \eqref{IERK: general form} if the two matrices $D_{\mathrm{E}}$ and $D_{\mathrm{EI}}$ are positive definite. Actually, it is easy to check that
\[
  \sum_{i=1}^k \sum_{j=1}^i\myinnerb{ d_{ij}(z) v^j, v^i} \geq \sum_{i=1}^k \sum_{j=1}^i \myinnerb{d_{ij}^{(e)} v^j, v^i}\quad\text{for $1\le k\le s$.}
\]
Thus the condition of Lemma \ref{lemma: bound quadratic form} can be fulfilled by taking $\lambda_{\es}$ as the minimum eigenvalue $\lambda_{\min}\kbrab{\mathcal{S}(D_{\mathrm{E}})}$ of the $z$-independent matrix $D_{\mathrm{E}}$, that is,
\begin{align}\label{IERK: choose lambda_es}     
	\lambda_{\es} := \lambda_{\min}\kbrab{\mathcal{S}(D_{\mathrm{E}})}.
\end{align}
 
By using the choice \eqref{IERK: choose lambda_es} of $\lambda_{\es}$, \cref{theorem: energy stability} is valid  for the IERK methods in \cite{FuTangYang:2024MCOM} and all of IERK methods developed recently in \cite{LiaoWangWen:2024IERK}, in which some parameterized  IERK methods are carefully designed to make the two matrices $D_{\mathrm{E}}$ and $D_{\mathrm{EI}}$ be positive definite. A preferred IERK method is then picked out by choosing the involved parameter such that the associated average dissipation rate (defined via the average eigenvalue of the differentiation matrix $D$) approaches 1 as close as possible.

\subsubsection{Second-order IERK methods}
Consider the one-parameter second-order IERK methods with the Butcher tableaux, see \cite[Eq. (3.6)]{LiaoWangWen:2024IERK},
\begin{equation}\label{scheme: IERK2}
	\begin{array}{c|ccc|ccc}
		0 & 0 &  && 0 &  &         \\
		\frac{\sqrt{2}}2 & \frac{\sqrt{2}}2-a_{33} & a_{33} & &   \frac{\sqrt{2}}2 & 0 &  \\[3pt]
		1 & \frac{\sqrt{2}-1+(2-\sqrt{2})a_{33}}{\sqrt{2}} & \frac{1-2a_{33}}{\sqrt{2}} & a_{33} & \frac{2-\sqrt{2}}2 &\frac{\sqrt{2}}2 & 0    \\[3pt]
		\hline  & \frac{\sqrt{2}-1+(2-\sqrt{2})a_{33}}{\sqrt{2}} & \frac{1-2a_{33}}{\sqrt{2}} & a_{33} & \frac{2-\sqrt{2}}2 &\frac{\sqrt{2}}2 & 0
	\end{array}.
\end{equation}
As proven in \cite[Theorem 3.2]{LiaoWangWen:2024IERK}, the  associated differentiation matrices of \eqref{scheme: IERK2} are positive definite if the parameter $a_{33}\ge\tfrac{1+\sqrt{2}}{4}$. 
Simple calculations yield that
\begin{equation*}
  D_{\mathrm{E}} = \begin{pmatrix}
   \sqrt{2} & 0 \\
   2 \sqrt{2}-2 & \sqrt{2}
  \end{pmatrix} \quad\text{and}\quad \lambda_{\min}\kbrab{\mathcal{S}(D_{\mathrm{E}})} = 1.
\end{equation*}
According to \eqref{IERK: choose lambda_es}, one can take $\lambda_{\es}:=1$ such that $\lambda_{\es} \leq \lambda_{\min}\kbrab{\mathcal{S}(D(z))}$ for $z<0$.

\subsubsection{Third-order IERK methods}
Consider the one-parameter third-order IERK methods with the Butcher tableaux, see \cite[Eq. (4.5)]{LiaoWangWen:2024IERK},
\begin{small}\begin{align}\label{Scheme: IERK3}
  \begin{array}{c|ccccc|ccccc}
     0 & 0 & & & & & 0 & \\[3pt]
    \tfrac{4}{5} & \tfrac{2}{25} & \tfrac{18}{25} & & & & \tfrac{4}{5} & 0  \\[3pt]
   \tfrac{7}{5} & \tfrac{3}{8} & \tfrac{61}{200} & \tfrac{18}{25} & & & \tfrac{3}{5} & \tfrac{4}{5} & 0   \\[3pt]
   \tfrac{6}{5} & a_{41} & a_{42} & a_{43} & \tfrac{18}{25} & & \tfrac{10111}{10080} & -\tfrac{6079}{10080} & \tfrac{4}{5} & 0  \\[3pt]
  1 & \tfrac{1030769}{2877000} & \tfrac{276523}{1233000} & -\tfrac{196127}{1078875} & -\tfrac{2068}{17125} & \tfrac{18}{25} & \tfrac{313}{840} & \tfrac{131}{360} & -\tfrac{169}{315} & \tfrac{4}{5} & 0 \\[3pt]
    \hline\\[-10pt]   & \tfrac{1030769}{2877000} & \tfrac{276523}{1233000} & -\tfrac{196127}{1078875} & -\tfrac{2068}{17125} & \tfrac{18}{25} & \tfrac{313}{840} & \tfrac{131}{360} & -\tfrac{169}{315} & \tfrac{4}{5} & 0
  \end{array}\,,
\end{align}
\end{small}
where $a_{41}=\tfrac{3 a_{43}}{4}+\tfrac{7277}{12600}$ and $\hat{a}_{42} = -\tfrac{7 a_{43}}{4}-\tfrac{1229}{12600}$. As shown in \cite[Theorem 4.2]{LiaoWangWen:2024IERK}, if $-0.633312 \leq a_{43} \leq -0.371114$, the associated differentiation matrices of \eqref{Scheme: IERK3} are positive definite. Simple calculations give  $\lambda_{\min}\kbrab{\mathcal{S}(D_{\mathrm{E}})} \approx 0.136355$ so that
one can choose $\lambda_{\es}:=0.13$ satisfying $\lambda_{\es} \leq \lambda_{\min}\kbrab{\mathcal{S}(D(z))}$ for $z<0$.


Actually, the mentioned Lobatto-type IERK methods \eqref{scheme: IERK2} and \eqref{Scheme: IERK3} are regarded as good candidates in \cite{LiaoWangWen:2024IERK} to preserve the original energy dissipation law for gradient flows. Similarly, the associated lower bounds $\lambda_{\es}$ for the parameterized Radau-type IERK methods in \cite{LiaoWangWen:2024IERK} can also be found and are left to interested readers. Note that, this delicate result \eqref{IERK: choose lambda_es}  comes from the fact that the differentiation matrix $D(z)$ is linear with respect to $z$. The situation becomes more complex for the EERK and CIFRK methods since the RK coefficients are always dependent on $z$.

\subsection{EERK methods}

When applying to the PFC type model \eqref{model: PFC type},  the EERK methods  \cite{HochbruckOstermann:2005SINUM,HochbruckOstermann:2010EERK} that naturally preserve the steady-state of \eqref{model: PFC type} take the form
\begin{align}\label{EERK: general form}     
	U^{n,i+1} = U^{n,1}-\tau\sum_{j=1}^{i}a_{i+1,j}({\tau}M_h L_{\kappa})M_h\kbrab{f_{\kappa}(U^{n,j})-L_{\kappa}U^{n,1}}
\end{align}
for $n\ge1$ and $1\le i\le s$, where 
the coefficients $a_{i+1,j}$ are constructed from linear combinations of the entire functions $\varphi_j(z)$ and scaled versions thereof. These functions are given by 
\begin{align*}
	\varphi_0(z)=e^z\quad\text{and}\quad\varphi_j(z):=\int_0^1e^{(1-s)z}\frac{s^{j-1}}{(j-1)!}\zd s \quad\text{for $z\in\mathbb{C}$ and $j\geq 1,$}
\end{align*}
which satisfy the recursion formula $\varphi_{k+1}(z)=\frac{\varphi_k(z)-1/k!}{z}$ for $k\geq 0$. The involved matrix functions $\varphi_j({\tau}M_hL_{\kappa})$ are defined on the spectrum of ${\tau}M_hL_{\kappa}$, that is, the values  $\{\varphi_j(\lambda_k): 1 \leq k \leq \mathrm{dim}({\tau}M_hL_{\kappa})\}$ exist, where $\lambda_k$ are the eigenvalues of ${\tau}M_hL_{\kappa}$ and thus $\varphi_j(\lambda_k)$ are the eigenvalues of $\varphi_j({\tau}M_hL_{\kappa})$. With the abscissas $c_1=0$ and $c_{s+1}=1$, in the following context, we also use the simplified notations $\varphi_{i,j}$,
\begin{align*}
	\varphi_{i,j}:=\varphi_{i,j}({\tau}M_hL_{\kappa})=\varphi_{i}(c_j{\tau}M_hL_{\kappa}),\quad i\ge0,\;1\leq j\leq s+1.
\end{align*}

As done before, we use the coefficient matrix $A(z):= \kbra{a_{i+1,j}(z)}_{i,j=1}^{s}$ to represent the EERK method \eqref{EERK: general form},
which can be formulated into the differential form \eqref{eq: unified differential form} with the associated differentiation matrix defined by, see \cite[Theorem 2.1]{FuShenYang:2024SCM,LiaoWang:2024MCOM},
\begin{align}\label{EERK: differentiation matrix}
	D(z) := A^{-1}(z) E_s + z E_s -\tfrac{z}2 I_s.
\end{align}
If the differentiation matrix $D(z)$ is positive definite for any $z<0$, \cref{theorem: energy stability} says that the corresponding EERK method preserves the original energy dissipation law in solving the PFC type model \eqref{model: PFC type}. Since the minimum eigenvalue of $\mathcal{S}(D(z))$ is dependent on $z$, the existence of a $z$-independent lower bound $\lambda_{\es}$ of the minimum eigenvalue $\lambda_{\min}\kbrab{\mathcal{S}(D(z))}$ should be further checked.


\subsubsection{Second-order EERK methods}

The parameterized EERK method \eqref{EERK: general form} satisfying the stiff-order two conditions has the following coefficient matrix, cf. \cite[Eq. (3.2)]{LiaoWang:2024MCOM} or \cite[Eq. (2.39)]{HochbruckOstermann:2010EERK}, 
\begin{equation}\label{scheme: EERK2}
	A(z; c_2) := \begin{pmatrix}
		c_{2}\varphi_{1,2} & 0\\[2pt]
		\varphi_1-\frac{1}{c_2}\varphi_2 & \frac{1}{c_2}\varphi_2 
	\end{pmatrix}\,.
\end{equation}
If the abscissa $c_2\in[\tfrac{1}{2}, 1]$, \cite[Corollary 3.1]{LiaoWang:2024MCOM} claims that the associated differentiation matrix $D(z; c_2)$ computed by \eqref{EERK: differentiation matrix} is positive definite for $z<0$. Here we consider the case $c_2=\tfrac12$ of \eqref{scheme: EERK2} since the associated average dissipation rate is closer to 1 than other cases $c_2\in(\tfrac{1}{2}, 1]$, see \cite[Figure 2, Table 1]{LiaoWang:2024MCOM}. 

\begin{figure}[htb!]
	\centering
	\subfigure[EERK2 method \eqref{scheme: EERK2} with $c_2=\tfrac12$]{\includegraphics[width=0.44\textwidth]{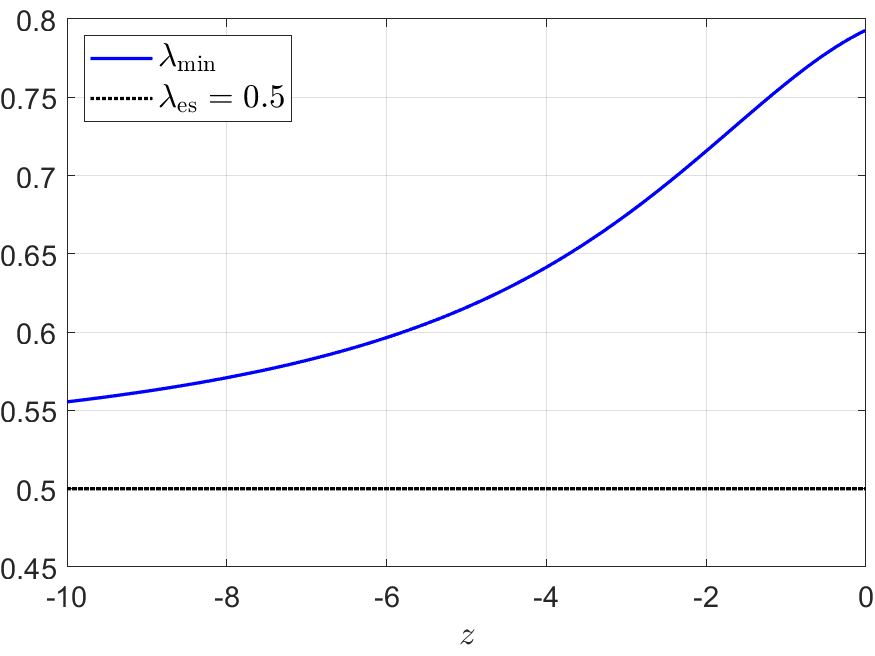}}
	\quad
	\subfigure[EERK2-w method \eqref{scheme: EERK2-w} with $c_2=\tfrac3{11}$]{\includegraphics[width=0.44\textwidth]{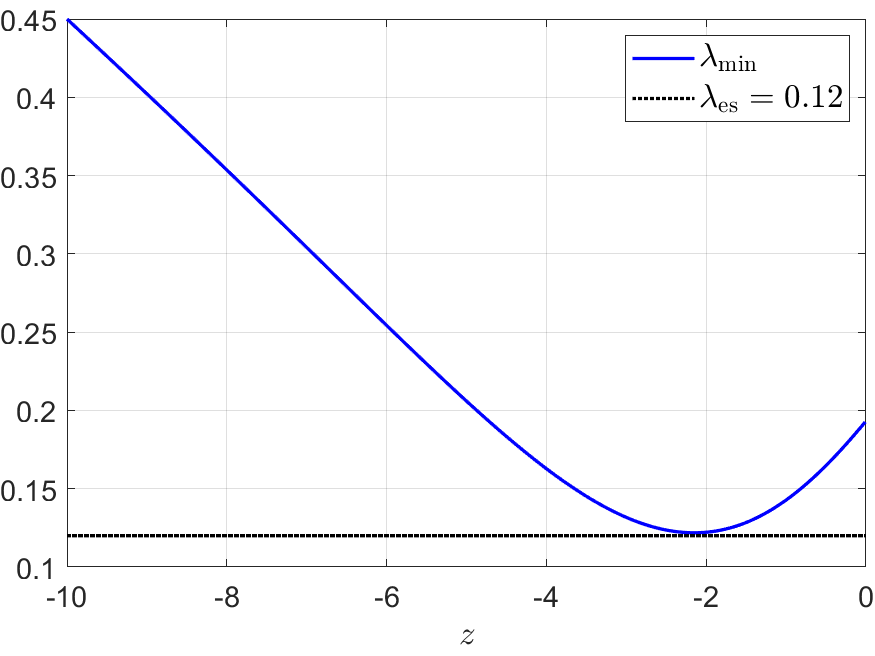}}
	\caption{Minimum eigenvalues of the second-order EERK methods.}
	\label{figure: minimum eigenvalue of EERK2}
\end{figure}

In general, a sharp bound of the minimum eigenvalue $\lambda_{\min}\kbrab{\mathcal{S}(D(z;\tfrac12))}$ would be mathematically complex and may be not essential to our current aim. As a compromise way, we will use a computer-aided approach, cf. \cite[Appendixes A and B]{LiaoWang:2024MCOM}, with the simple fact $\lim_{z\rightarrow-\infty}z^{p_1}e^{p_2z}=0$ for two constants $p_2>0$ and $p_1\ge0$.  
Actually, it is not difficult to  find that the minimum eigenvalue has the asymptotic behavior, $\lambda_{\min}\kbrab{\mathcal{S}(D(z;\tfrac12))}\sim\tfrac{z}{2 z+2}$ as $z\to-\infty$; while \cref{figure: minimum eigenvalue of EERK2}(a) says that the curve $\lambda_{\min}\kbrab{\mathcal{S}(D(z;\tfrac12))}$ is increasing in $z\in[-10,0]$. Thus, one can choose $\lambda_{\es} = \tfrac12$ such that $\lambda_{\es} \le\lambda_{\min}\kbrab{\mathcal{S}(D(z;\tfrac12))}$ for $z<0$. It is easy to check that, 
this lower bound $\lambda_{\es} = \tfrac12$ is also valid for \lan{case} $c_2=1$ of \eqref{scheme: EERK2},
which corresponds to the ETDRK2 method discussed in \cite{LiQiaoWangZheng:2024arxiv} for the PFC equation.

As the weak version  (stiff-order one) of \eqref{scheme: EERK2}, second-order EERK-w method has the following coefficient matrix, cf. \cite[Eq. (3.5)]{LiaoWang:2024MCOM} or \cite[Eq. (2.40)]{HochbruckOstermann:2010EERK},
\begin{equation}\label{scheme: EERK2-w}
	A^{(w)}(z;c_2) := \begin{pmatrix}
		c_{2}\varphi_{1,2} & 0\\[2pt]
		\varphi_1-\frac{1}{2c_2}\varphi_1 & \frac{1}{2c_2}\varphi_1
	\end{pmatrix}\,.
\end{equation}
The associated differentiation matrix $D^{(w)}(z; c_2)$ of the EERK-w methods \eqref{scheme: EERK2-w} with the abscissa $c_2\in[\tfrac{3}{11}, 1]$ was shown in \cite[Corollary 3.2]{LiaoWang:2024MCOM} to be positive definite for any $z<0$. As illustrated in \cref{figure: minimum eigenvalue of EERK2}(b), the curve of the minimum eigenvalue $\lambda_{\min}\kbrab{\mathcal{S}(D^{(w)}(z;\tfrac3{11}))}$ is decreasing for $z \leq -4$, while the associated asymptotic expression is $\lambda_{\min}\kbrab{\mathcal{S}(D^{(w)}(z;\tfrac3{11}))}\sim-\tfrac{z}{22}$ when $z\to-\infty$. By simple but rather lengthy calculations, one can verify that $\lambda_{\min}\kbrab{\mathcal{S}(D^{(w)}(z;\tfrac3{11}))}$ achieves the minimum value $0.121853$ at $z=-2.15036$. Thus one can choose $\lambda_{\es} = 0.12$ so that $\lambda_{\es} \le \lambda_{\min}\kbrab{\mathcal{S}(D^{(w)}(z;\tfrac3{11}))}$ for all $z<0$.

The other cases $c_2\in(\tfrac{3}{11}, 1]$ of the EERK-w methods \eqref{scheme: EERK2-w} can be handled similarly. Typically, for the case $c_2=\tfrac12$ discussed in \cite{SunZhangetal:2024arxiv} for the SH equation,  one can choose $\lambda_{\es} = 0.79$ so that $\lambda_{\es} \le \lambda_{\min}\kbrab{\mathcal{S}(D^{(w)}(z;\tfrac1{2}))}$ for all $z<0$. 



\subsubsection{Third-order EERK methods}

We consider the one-parameter family of third-order EERK3-1 method, see \cite[Eq. (4.2)]{LiaoWang:2024MCOM} or \cite[Eq. (5.8)]{HochbruckOstermann:2005SINUM},
\begin{align}\label{scheme: EERK3-1}
	A(z; c_2) = \begin{pmatrix}
		c_{2}\varphi_{1,2} &   \\[2pt]
		\frac{2}{3}\varphi_{1,3}-\frac{4}{9c_2}\varphi_{2,3} &  \frac{4}{9c_2}\varphi_{2,3} \\[2pt]
		\varphi_1-\frac{3}{2}\varphi_2 & 0 & \frac{3}{2}\varphi_2
	\end{pmatrix}\,,
\end{align}
and the two-parameter EERK3-2 method, see \cite[Eq. (4.3)]{LiaoWang:2024MCOM} or \cite[Eq. (5.9)]{HochbruckOstermann:2005SINUM},
\begin{align}\label{scheme: EERK3-2}
	A(z; c_2, c_3) = \begin{pmatrix}
		c_{2}\varphi_{1,2} & \\[2pt]
		c_{3}\varphi_{1,3}-a_{32} &  \gamma c_2\varphi_{2,2}+\frac{c_3^2}{c_2}\varphi_{2,3} \\[2pt]
		\varphi_1-a_{42}-a_{43} & \frac{\gamma}{\gamma c_2+c_3}\varphi_2 & \frac{1}{\gamma c_2+c_3}\varphi_2
	\end{pmatrix}\,,
\end{align}
where the parameter $\gamma:=\frac{(3c_3-2)c_3}{(2-3c_2)c_2}$ for $c_2\neq\frac{2}{3}$ and $c_2\neq c_3$ (to ensure $a_{32}\neq0$). 
For the two third-order EERK methods, the corresponding differentiation matrices become more complex than those of second-order EERK methods. 

\begin{figure}[htb!]
	\centering
	\subfigure[EERK3-1 method \eqref{scheme: EERK3-1} with $c_2=\tfrac{4}{9}$]{\includegraphics[width=0.44\textwidth]{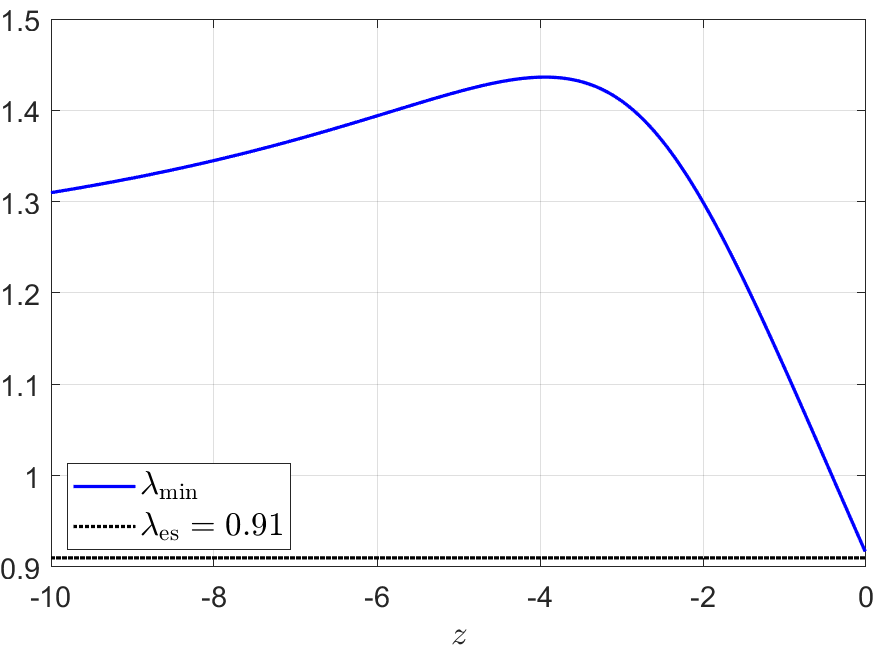}}
	\quad
	\subfigure[EERK3-2 \eqref{scheme: EERK3-2} with $(c_2, c_3)=(\tfrac{1}2, \tfrac{7}{10})$]{\includegraphics[width=0.44\textwidth]{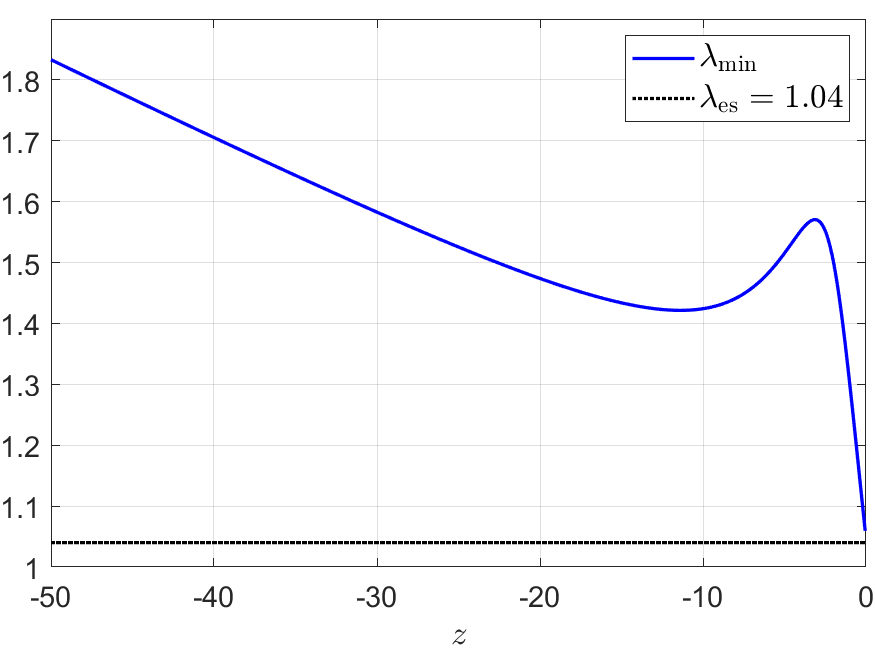}}
	\caption{Minimum eigenvalues of the third-order EERK methods.}
	\label{figure: minimum eigenvalue of EERK3}
\end{figure}

If the abscissa $c_2\in[\tfrac{4}{9}, 1]$, the differentiation matrices of the EERK3-1 methods \eqref{scheme: EERK3-1} are positive definite for $z<0$, see \cite[Corollary 4.1]{LiaoWang:2024MCOM}. Here we consider the case $c_2=\tfrac{4}9$ of \eqref{scheme: EERK3-1} since the associated average dissipation rate is closer to 1 than other cases $c_2\in(\tfrac{4}{9}, 1]$, cf. \cite[Figure 5, Table 2]{LiaoWang:2024MCOM}. \cref{figure: minimum eigenvalue of EERK3}(a) depicts the minimum eigenvalue of the associated differentiation matrix $D(z;\tfrac{4}9)$ of the EERK3-1 method \eqref{scheme: EERK3-1} for the case $c_2=\tfrac{4}{9}$. Simple calculations yield that $ \lambda_{\min}\kbrab{\mathcal{S}(D(z;\tfrac{4}9))}\approx 0.914971$ as $z\to0$, while  $\lambda_{\min}\kbrab{\mathcal{S}(D(z;\tfrac{4}9))}\sim\tfrac{5}{12}-\tfrac{z+\sqrt{9 z+z^{2}}}{6}$ as $z\to-\infty$. Thus, one can choose $\lambda_{\es}:=0.91$ so that $\lambda_{\es}\le\lambda_{\min}\kbrab{\mathcal{S}(D(z;\tfrac{4}9))}$ for $z<0$. 


In \cite[Corollaries 4.2-4.4]{LiaoWang:2024MCOM}, the  differentiation matrices of the EERK3-2 \lan{methods} \eqref{scheme: EERK3-2} for three special cases, $(c_2, c_3)=(1, \tfrac{1}{2})$, $(c_2, c_3)=(\tfrac{3}4, \tfrac{3}5)$ and $(c_2, c_3)=(\tfrac{1}{2}, \tfrac{7}{10})$, are proven to be positive definite for  $z<0$. 
\cref{figure: minimum eigenvalue of EERK3}(b) depicts the minimum eigenvalue of the associated differentiation \lan{matrix} $D(z;\tfrac{1}2, \tfrac{7}{10})$ of the case $(c_2, c_3)=(\tfrac{1}{2}, \tfrac{7}{10})$, which has the smallest value of average dissipation rate among these three cases, cf. \cite[Table 2]{LiaoWang:2024MCOM}. One has $\lambda_{\min}\kbrab{\mathcal{S}(D(z;\tfrac{1}2, \tfrac{7}{10}))}\approx 1.04656$ as $z\to0$; and for proper large $\abst{z}$, the minimum eigenvalue has the asymptotic expression,
$$\lambda_{\min}\kbrab{\mathcal{S}(D(z;\tfrac{1}2, \tfrac{7}{10}))}\sim-\tfrac{\sqrt{183698}-457}{2100}z \approx -0.0135238z\quad\text{as $z\to-\infty$}.$$
We see that the minimum eigenvalue attains the smallest value at $z=0$. Then one can choose $\lambda_{\es}:=1.04$ so that $\lambda_{\es}\le\lambda_{\min}\kbrab{\mathcal{S}(D(z;\tfrac{1}2, \tfrac{7}{10}))}$ for any $z<0$.


\subsection{CIFRK methods}
The CIFRK methods are suggested recently in \cite{LiaoWangWen:2024JCP} to remedy the so-called exponential effect and the non-preservation of original energy dissipation of the IFRK (known as Lawson \cite{Lawson:1967SINUM}) methods, when applying to \eqref{model: PFC type},
\begin{align}\label{scheme: Lawson}
	u_h^{n,i+1}= e^{c_{i+1} \tau M_h L_{\kappa}} u_h^{n,1}-\tau \sum_{j=1}^{i}\hat{a}_{i+1,j}(\tau M_hL_{\kappa})f_{\kappa}(u_h^{n, j})
\end{align}
for $n\ge1$ and $ 1\le i\le s$. The RK coefficients $\hat{a}_{i+1,j}(z):=\hat{a}_{i+1,j} (0)e^{(c_{i+1}-c_{j})z},$ where $\hat{a}_{i+1,j}(0)$ are real numbers that are independent of any discrete operators including $M_h$ and $L_\kappa$. The IFRK method \eqref{scheme: Lawson} can be represented by the lower triangular matrices  $A_{\mathrm{E}}(z):=\kbrab{\hat{a}_{i+1,j}(z)}_{i,j=1}^{s}$; while $A_{\mathrm{E}}(0):=\kbrab{\hat{a}_{i+1,j}(0)}_{i,j=1}^{s}$ is the coefficient matrix of the underlying explicit RK method.
	
By enforcing the preservation of steady-state for any time-step sizes, the telescopic correction (by modifying the coefficient of $u_h^{n,i+1}$) and the nonlinear-term translation correction (by modifying the coefficient of $f_{\kappa}(u_h^{n, i})$) applied to a given IFRK method \eqref{scheme: Lawson} arrive at the TIF and NIF methods, respectively, see \cite[Section 3]{LiaoWangWen:2024JCP} for more details. The two classes of CIFRK methods can be formulated as the following steady-state preserving forms  
\begin{align}
  \text{TIF:}\quad& u_h^{n,i+1}= u_h^{n,1} - \tau\sum_{j=1}^i\hat{a}^{\cte}_{i+1,j}(\tau M_h L_{\kappa}) M_h \kbrab{f_{\kappa}(u_h^{n, j}) - L_{\kappa} u_h^{n,1}}, \label{scheme: TIF} \\
  \text{NIF:}\quad& u_h^{n,i+1}= u_h^{n,1} - \tau\sum_{j=1}^i\hat{a}^{\cnt}_{i+1,j}(\tau M_h L_{\kappa}) M_h \kbrab{f_{\kappa}(u_h^{n, j}) - L_{\kappa} u_h^{n,1}} \label{scheme: NIF}
\end{align}
for $n\ge1$ and $1\le i\le s$. We will use the coefficient matrices $A_{\cte}(z):=\kbrab{\hat{a}_{i+1,j}^{\cte}(z)}_{i,j=1}^{s}$ and $A_{\cnt}(z):=\kbrab{\hat{a}_{i+1,j}^{\cnt}(z)}_{i,j=1}^{s}$ to \lan{represent} the TIF and NIF methods, respectively. Then one can get the associated differential form \eqref{eq: unified differential form} with the following differentiation matrices for the two classes of CIFRK methods, respectively,
\begin{align}\label{CIFRK: differentiation matrix}
  &\,D_{\cte}(z) := A_{\cte}^{-1}(z) E_s + z E_s -\tfrac{z}2 I_s \quad\text{and}\quad
 D_{\cnt}(z) := A_{\cnt}^{-1}(z) E_s + z E_s -\tfrac{z}2 I_s.
  \end{align}

It is easy to know that \cref{theorem: energy stability} holds for the CIFRK methods \eqref{scheme: TIF}-\eqref{scheme: NIF} if the matrices $D_{\cte}(z)$ and $D_{\cnt}(z)$ in \eqref{CIFRK: differentiation matrix} are positive definite for any $z<0$. In general, the minimum eigenvalues of $\mathcal{S}(D_{\cte}(z))$ and $\mathcal{S}(D_{\cnt}(z))$ are dependent on $z$. Thus the existence of a $z$-independent lower bound $\lambda_{\es}$ of these minimum eigenvalues should be further checked. In the following, eight CIFRK schemes stemmed from the Heun's and Ralston's explicit RK schemes \cite{Heun:1900ZMP,Ralston:1962} will be revisited.



\subsubsection{Second-order CIFRK methods}

The coefficient matrices $A_{\mathrm{E}}^{(2,H)}(0)$ and $A_{\mathrm{E}}^{(2,R)}(0)$ of the Heun's and Ralston's explicit methods read, respectively,
\begin{align*}
	A_{\mathrm{E}}^{(2,H)}(0) := \begin{pmatrix}
		1 & 0 \\[2pt]
		\tfrac{1}{2} & \tfrac{1}{2}
	\end{pmatrix}\,\quad\text{and}\quad
	 A_{\mathrm{E}}^{(2,R)}(0) := \begin{pmatrix}
		\tfrac{2}{3} & 0 \\[2pt]
		\tfrac{1}{4}& \tfrac{3}{4}
	\end{pmatrix}\,. 
\end{align*}

For the second-order Heun scheme, the resulting TIF2-Heun and NIF2-Heun methods have the following coefficient matrices \cite[Eq. (3.11) \& (3.13)]{LiaoWangWen:2024JCP}
\begin{align}
  \text{TIF2-Heun:}\quad& A_{\cte}^{(2,H)}(z) := \begin{pmatrix}
    \tfrac{1}{1-z} & 0 \\[2pt]
    \tfrac{1}{2 -z(1+e^{-z})} & \tfrac{1}{2e^{z} -z(1+e^{z})}
  \end{pmatrix}\,, \label{TIF2: Heun} \\
  \text{NIF2-Heun:}\quad& A_{\cnt}^{(2,H)}(z) := \begin{pmatrix}
    \tfrac{e^{z}-1}{z} & 0 \\[2pt]
    \tfrac{1}{2}e^{z} & \tfrac{e^{z}-1}{z}-\tfrac{1}{2}e^{z}
  \end{pmatrix}\,; \label{NIF2: Heun}
\end{align}
while the associated differentiation matrices $D_{\cte}^{(2,H)}(z)$ and $D_{\cnt}^{(2,H)}(z)$ can be computed via the formulas in \eqref{CIFRK: differentiation matrix}. 

\begin{figure}[htb!]
  \centering
  \subfigure[CIF2-Heun]{\includegraphics[width=0.44\textwidth]{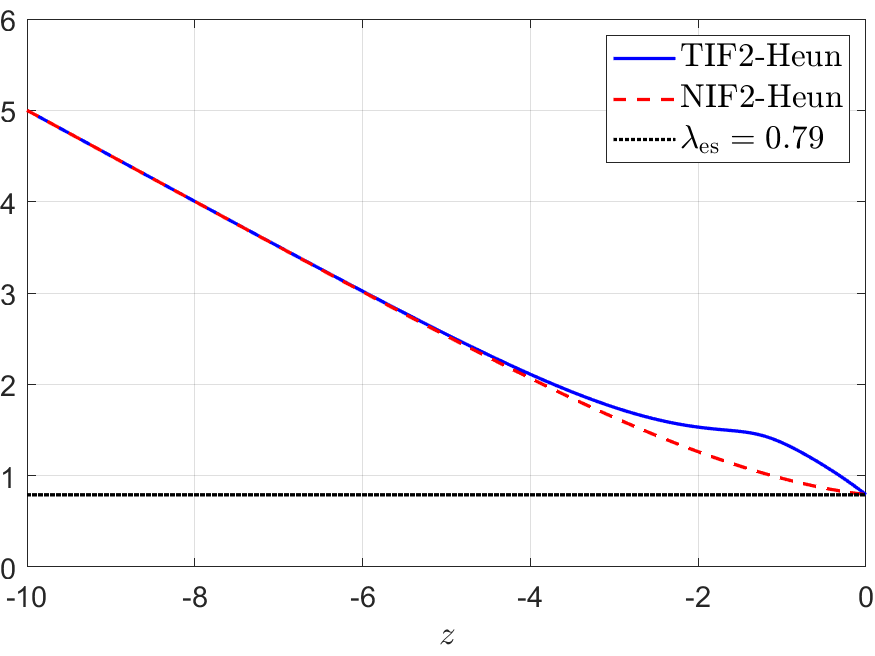}}
  \quad
  \subfigure[CIF2-Ralston]{\includegraphics[width=0.44\textwidth]{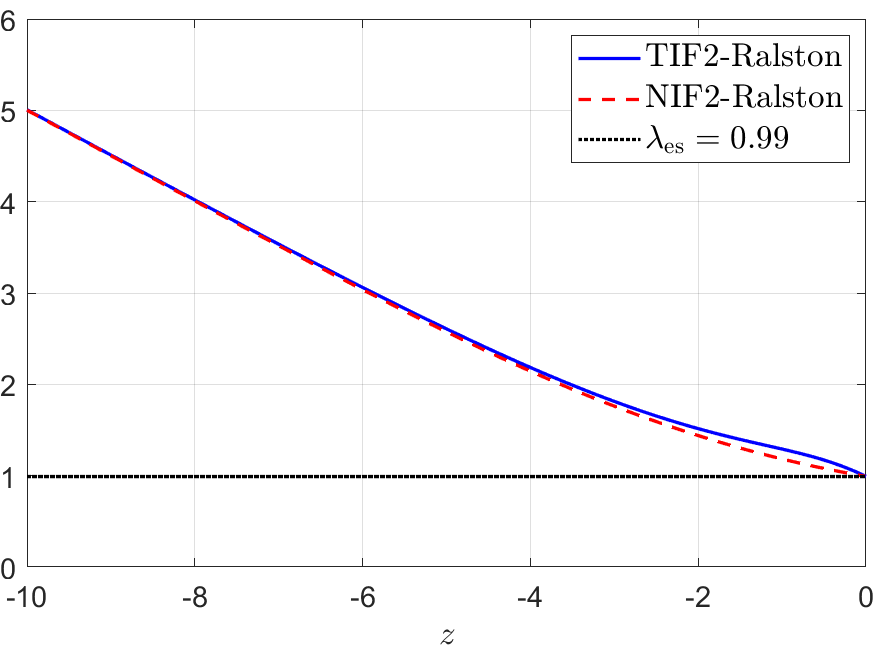}}
  \caption{Minimum eigenvalues of the second-order CIFRK schemes.}
  \label{figure: minimum eigenvalue of CIF2}
\end{figure}

In \cite[Theorem 3.1]{LiaoWangWen:2024JCP}, the differentiation matrices $D_{\cte}^{(2,H)}(z)$ and $D_{\cnt}^{(2,H)}(z)$ are shown to be positive definite for $z<0$. 
It is easy to check that 
$$\lambda_{\min}\kbrab{\mathcal{S}(D_{\cte}^{(2,H)})}\sim\lambda_{\min}\kbrab{\mathcal{S}(D_{\cnt}^{(2,H)})}\sim-\tfrac{z}2\quad\text{as $z\to-\infty$}.$$
They indicate the decreasing behavior of the minimum eigenvalue $\lambda_{\min}=\lambda_{\min}(z)$ when $|z|$ is large.
In \cref{figure: minimum eigenvalue of CIF2}(a), we depict the curves of $\lambda_{\min}\kbrab{\mathcal{S}(D_{\cte}^{(2,H)})}$ and $\lambda_{\min}\kbrab{\mathcal{S}(D_{\cnt}^{(2,H)})}$  near $z=0$. As seen, the two curves of minimum eigenvalues are decreasing with respect to $z$.  Thus one can choose $\lambda_{\es}:= 0.79$ such that
 \begin{align*}
 \lambda_{\es}\le\min\big\{\lambda_{\min}\kbrab{\mathcal{S}(D_{\cte}^{(2,H)})},\lambda_{\min}\kbrab{\mathcal{S}(D_{\cnt}^{(2,H)})}\big\}\quad\text{for $z<0$,}
 \end{align*}
while the latter approaches $\tfrac{3-\sqrt{2}}{2}\approx 0.792893$ as $z \to 0$.

For the second-order Ralston scheme, the resulting TIF2-Ralston and NIF2-Ralston methods have the following coefficient matrices \cite[Eq. (3.15) \& (3.17)]{LiaoWangWen:2024JCP}
\begin{align}
  \text{TIF2-Ralston:}\quad& A_{\cte}^{(2,R)} = \begin{pmatrix}
    \tfrac{2}{3-2z} & 0 \\[2pt]
    \tfrac{e^{\frac{2z}{3}}}{(4-z)e^{\frac{2z}{3}}-3z} & \tfrac{3}{(4-z)e^{\frac{2z}{3}}-3z}
  \end{pmatrix}\,, \label{TIF2: Ralston} \\
  \text{NIF2-Ralston:}\quad& A_{\cnt}^{(2,R)} = \begin{pmatrix}
    \tfrac{e^{\frac{2 z}{3}}-1}{z} & 0 \\[2pt]
    \tfrac{1}{4}e^z & \tfrac{e^z-1}{z}-\tfrac{1}{4}e^z
  \end{pmatrix}\,. \label{NIF2: Ralston}
\end{align}
In \cite[Theorem 3.2]{LiaoWangWen:2024JCP}, the associated differentiation matrices $D_{\cte}^{(2,R)}(z)$ and $D_{\cnt}^{(2,R)}(z)$ are proven to be positive definite for $z<0$.
The minimum eigenvalues of $\mathcal{S}(D_{\cte}^{(2,R)})$ and $\mathcal{S}(D_{\cnt}^{(2,R)})$ also have the same asymptotic behavior, $\lambda_{\min}\sim-\tfrac{z}{2}$ when $|z|$ is large, while as depicted in \cref{figure: minimum eigenvalue of CIF2}(b), they are always decreasing as $z\in[-10, 0]$. Thus, one can choose $\lambda_{\es}:= 0.99$ such that
\begin{align*}
	\lambda_{\es}\le\min\big\{\lambda_{\min}\kbrab{\mathcal{S}(D_{\cte}^{(2,R)})},\lambda_{\min}\kbrab{\mathcal{S}(D_{\cnt}^{(2,R)})}\big\}\quad\text{for $z<0$,}
\end{align*}
while the latter approaches the value $\tfrac{17-\sqrt{26}}{12}\approx 0.991748$ as $z \to 0$.

\subsubsection{Third-order CIFRK methods} 
The coefficient matrices $A_{\mathrm{E}}^{(3,H)}(0)$ and $A_{\mathrm{E}}^{(3,R)}(0)$ of the third-order Heun's and Ralston's RK methods read, respectively,
\begin{align*}
	A_{\mathrm{E}}^{(3,H)}(0) := \begin{pmatrix}
		\frac{1}{3}& & \\[1pt]
		0 & \frac{2}{3}&\\[1pt]
		\frac{1}{4} & 0 & \frac{3}{4}
	\end{pmatrix}\,\quad\text{and}\quad
	A_{\mathrm{E}}^{(3,R)}(0) := \begin{pmatrix}
		\frac{1}{2} & & \\[1pt]
		0 & \frac{3}{4}&\\[1pt]
		\frac{2}{9} & \frac{1}{3} & \frac{4}{9}
	\end{pmatrix}\,. 
\end{align*}

\begin{figure}[htb!]
	\centering
	\subfigure[CIF3-Heun]{\includegraphics[width=0.44\textwidth]{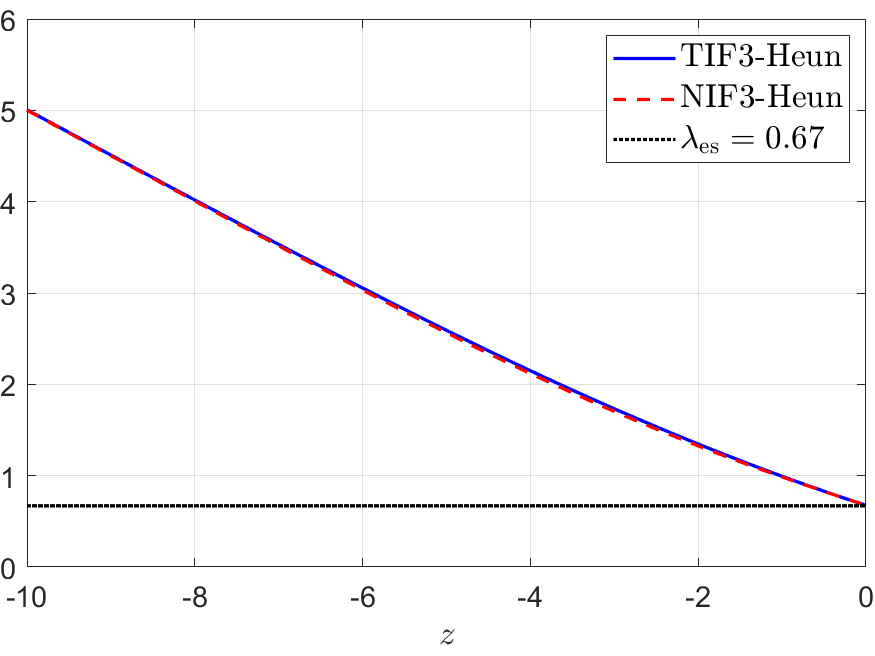}}
	\quad
	\subfigure[CIF3-Ralston]{\includegraphics[width=0.44\textwidth]{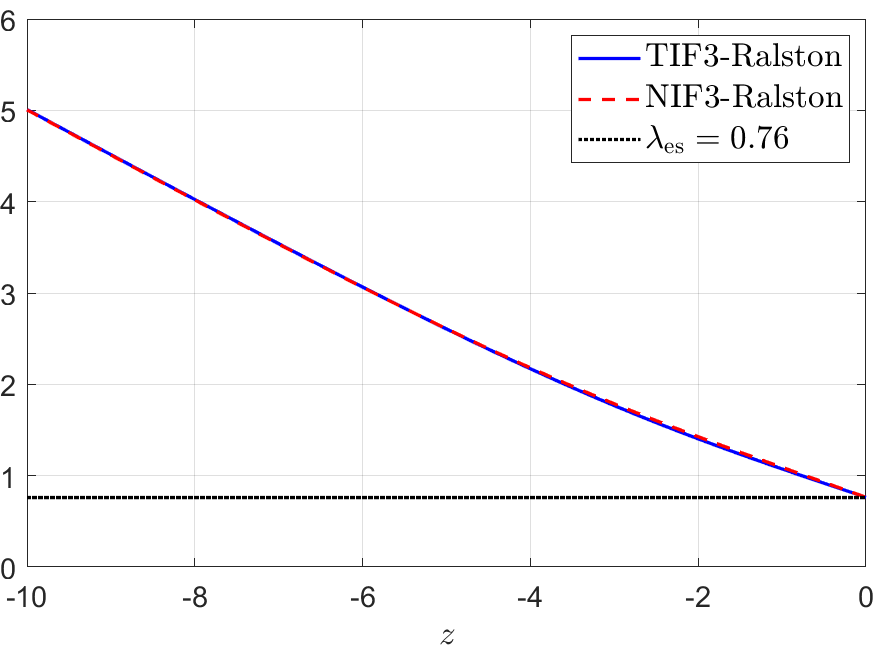}}
	\caption{Minimum eigenvalues of the third-order CIFRK methods.}
	\label{figure: minimum eigenvalue of CIFRK3}
\end{figure}

For the third-order Heun's explicit RK scheme, the resulting TIF3-Heun and NIF3-Heun methods have the following coefficient matrices \cite[Eqs. (4.1)-(4.2)]{LiaoWangWen:2024JCP}
\begin{align*}
  & A_{\cte}^{(3,H)}(z) := \begin{pmatrix}
   \frac{1}{3-z} \\[2pt]
   0 & \frac{2}{3 e^{\frac{z}{3}}-2 z}  \\[2pt]
   \frac{e^{\frac{2 z}{3}}}{e^{\frac{2 z}{3}} (4-z)-3 z} & 0 & \frac{3}{e^{\frac{2 z}{3}} (4-z)-3 z}
  \end{pmatrix}\,, \\
  & A_{\cnt}^{(3,H)}(z):= \begin{pmatrix}
   \frac{e^{\frac{z}{3}}-1}{z} \\[2pt]
   0 & \frac{e^{\frac{2 z}{3}}-1}{z} \\[2pt]
   \frac{e^z}{4} & 0 & \frac{e^z-1}{z}-\frac{e^z}{4}
  \end{pmatrix}\,.
\end{align*}
In \cite[Theorem 4.1]{LiaoWangWen:2024JCP}, the associated differentiation matrices $D_{\cte}^{(3,H)}(z)$ and $D_{\cnt}^{(3,H)}(z)$ are shown to be positive definite for $z<0$.
Also,  the minimum eigenvalues of $\mathcal{S}(D_{\cte}^{(3,H)})$ and $\mathcal{S}(D_{\cnt}^{(3,H)})$ have the same asymptotic behavior, $\lambda_{\min}\sim-\tfrac{z}{2}$ as $z\to-\infty$. \cref{figure: minimum eigenvalue of CIFRK3}(a) says that the curves $\lambda_{\min}\kbrab{\mathcal{S}(D_{\cte}^{(3,H)})}$ and $\lambda_{\min}\kbrab{\mathcal{S}(D_{\cnt}^{(3,H)})}$ are decreasing for $z\in[-10, 0]$. Then we can choose $\lambda_{\es}:= 0.67$ such that
\begin{align*}
	\lambda_{\es}\le\min\big\{\lambda_{\min}\kbrab{\mathcal{S}(D_{\cte}^{(3,H)})},\lambda_{\min}\kbrab{\mathcal{S}(D_{\cnt}^{(3,H)})}\big\}\quad\text{for $z<0$,}
\end{align*}
since the latter approaches $0.675972$ as $z\to 0$.

For the third-order Ralston's RK scheme, the resulting  TIF3-Ralston and NIF3-Ralston methods have the following coefficient matrices \cite[Eqs. (4.3)-(4.4)]{LiaoWangWen:2024JCP}
\begin{align*}
	&A_{\cte}^{(3,R)}(z) := \begin{pmatrix}
		\frac{1}{2-z}  \\[2pt]
		0 & \frac{3}{4 e^{\frac{z}{2}}-3 z}  \\[2pt]
		\frac{2 e^{\frac{3 z}{4}}}{e^{\frac{3 z}{4}} (9-2 z)-3 e^{\frac{z}{4}} z-4 z} & \frac{3 e^{\frac{z}{4}}}{e^{\frac{3 z}{4}} (9-2 z)-3 e^{\frac{z}{4}} z-4 z} & \frac{4}{e^{\frac{3 z}{4}} (9-2 z)-3 e^{\frac{z}{4}} z-4 z}
	\end{pmatrix}\,, \\
	&A_{\cnt}^{(3,R)}(z):= \begin{pmatrix}
		\frac{e^{\frac{z}{2}}-1}{z}  \\[2pt]
		0 & \frac{e^{\frac{3 z}{4}}-1}{z}  \\[2pt]
		\frac{2 e^z}{9} & \frac{e^{\frac{z}{2}}}{3} & \frac{e^z-1}{z}-\frac{e^{\frac{z}{2}}}{3}-\frac{2 e^z}{9}
	\end{pmatrix}\,.
\end{align*}
As stated by \cite[Theorem 4.1]{LiaoWangWen:2024JCP}, the associated differentiation matrices $D_{\cte}^{(3,R)}(z)$ and $D_{\cnt}^{(3,R)}(z)$ are positive definite for $z<0$. It is checked that the minimum eigenvalues of $\mathcal{S}(D_{\cte}^{(3,R)})$ and $\mathcal{S}(D_{\cnt}^{(3,R)})$ have the same asymptotic formula, $\lambda_{\min}\sim-\tfrac{z}{2}$ as $z\to-\infty$. \cref{figure: minimum eigenvalue of CIFRK3}(b) says that the curves $\lambda_{\min}\kbrab{\mathcal{S}(D_{\cte}^{(3,R)})}$ and $\lambda_{\min}\kbrab{\mathcal{S}(D_{\cnt}^{(3,R)})}$ are decreasing for $z\in[-10, 0]$. We choose $\lambda_{\es}:= 0.76$ such that
\begin{align*}
	\lambda_{\es}\le\min\big\{\lambda_{\min}\kbrab{\mathcal{S}(D_{\cte}^{(3,R)})},\lambda_{\min}\kbrab{\mathcal{S}(D_{\cnt}^{(3,R)})}\big\}\quad\text{for $z<0$,}
\end{align*}
while the latter approaches $0.763096$ as $z\to 0$.

As the end of this subsection, we claim that \cref{theorem: energy stability} is always valid for the eight CIFRK methods exhibited in \cite{LiaoWangWen:2024JCP} up to third-order accuracy.
It is interesting to note that the elements of the coefficient matrices $A_{\cte}(z)$ and $A_{\cnt}(z)$ of the above eight CIFRK methods are nonnegative for any $z<0$, cf. \cite[Figures 5, 6, 10 and 11]{LiaoWangWen:2024JCP}, which naturally inherit the nonnegative property of RK coefficients in the underlying Heun's and Ralston's explicit RK schemes. 

\section{Conclusion}\label{section: conclusion}

We present a concise theoretical framework for the original energy dissipation \lan{laws} of the IERK, EERK and CIFRK methods applying to the SH and PFC models. As the  main result, Theorem \ref{theorem: energy stability} establishes the global-in-time energy estimate and global-in-time uniform boundedness of stage solutions for three efficient classes of  RK methods for the PFC type model \eqref{model: PFC type}, under a simple condition that the associated  differentiation matrix $D(z)$ is positive definite for $z<0$.

In addition to the positive definiteness condition for the associated differentiation matrix of RK methods, the unified framework in \cref{section: theoretical framework} is stemmed from the two physical structures of the underlying phase field model: (PHY-1) the boundedness of discrete energy functional $E[u^{n,i}]$ implies the $H^2$ norm boundedness of the solution $u_h^{n,i}$; (PHY-2) the nonlinear bulk $f$ is a smooth function of the phase variable $\Phi$. However, these structures would be not fulfilled for other phase field models:
\begin{itemize}[leftmargin=10mm]
	\item[(a)] The Cahn-Hilliard model \cite{ChengWangWiseYue:2016JSC,DongLiQiaoZhang:2023ANM} with the Ginzburg-Landau free energy satisfy (PHY-2), but do not satisfy (PHY-1) since the boundedness of discrete energy functional implies only the  $H^1$ norm boundedness of the solution.
	\item[(b)] The molecular beam epitaxy equation with slope selection \cite{MoldovanGolubovic:2000MBE} and the square phase field crystal equation \cite{ElderGrant:2004SPFC} satisfy (PHY-1) but do not satisfy (PHY-2) since the involved nonlinear bulk $f=f(\nabla \Phi)$ is a function of the gradient.
\end{itemize}
For the nonlinear phase field models who cannot simultaneously satisfy the two physical properties (PHY-1) and (PHY-2), the $L^2$ norm and $H^1$ semi-norm boundedness in Lemma \ref{lemma: Energy implies H2 norm} for the nonlinear bulk $f$ will be no longer valid. Thus, whether a unified theory on the original energy dissipation law can be established for these nonlinear equations remains to be further explored.

\section*{Acknowledgment}
The authors would like to thank Dr. Bingquan Ji for his sincere discussions on \cref{lemma: Energy implies H2 norm}.

\bibliographystyle{siamplain}
\bibliography{referencesAdd20241210}

\end{document}